
%
\magnification1200
\pretolerance=100
\tolerance=200
\hbadness=1000
\vbadness=1000
\linepenalty=10
\hyphenpenalty=50
\exhyphenpenalty=50
\binoppenalty=700
\relpenalty=500
\clubpenalty=5000
\widowpenalty=5000
\displaywidowpenalty=50
\brokenpenalty=100
\predisplaypenalty=7000
\postdisplaypenalty=0
\interlinepenalty=10
\doublehyphendemerits=10000
\finalhyphendemerits=10000
\adjdemerits=160000
\uchyph=1
\delimiterfactor=901
\hfuzz=0.1pt
\vfuzz=0.1pt
\overfullrule=5pt
\hsize=146 true mm
\vsize=8.9 true in
\maxdepth=4pt
\delimitershortfall=.5pt
\nulldelimiterspace=1.2pt
\scriptspace=.5pt
\normallineskiplimit=.5pt
\mathsurround=0pt
\parindent=20pt
\catcode`\_=11
\catcode`\_=8
\normalbaselineskip=12pt
\normallineskip=1pt plus .5 pt minus .5 pt
\parskip=6pt plus 3pt minus 3pt
\abovedisplayskip = 12pt plus 5pt minus 5pt
\abovedisplayshortskip = 1pt plus 4pt
\belowdisplayskip = 12pt plus 5pt minus 5pt
\belowdisplayshortskip = 7pt plus 5pt
\normalbaselines
\smallskipamount=\parskip
 \medskipamount=2\parskip
 \bigskipamount=3\parskip
\jot=3pt
%
%
\def\ref#1{\par\noindent\hangindent2\parindent
 \hbox to 2\parindent{#1\hfil}\ignorespaces}
%
%
\font\tenss=cmss10
\font\sevenss=cmss8 at 7pt
\font\fivess=cmss8 at 5pt
\newfam\ssfam %
\textfont\ssfam=\tenss
\scriptfont\ssfam=\sevenss
\scriptscriptfont\ssfam=\fivess
%
%
%
%
%
%
%
%
%
\catcode`\_=11
\def\suf_fix{}
\def\scaled_rm_box#1{%
 \relax
 \ifmmode
   \mathchoice
    {\hbox{\tenrm #1}}%
    {\hbox{\tenrm #1}}%
    {\hbox{\sevenrm #1}}%
    {\hbox{\fiverm #1}}%
 \else
  \hbox{\tenrm #1}%
 \fi}
\def\suf_fix_def#1#2{\expandafter\def\csname#1\suf_fix\endcsname{#2}}
\def\I_Buchstabe#1#2#3{%
 \suf_fix_def{#1}{\scaled_rm_box{I\hskip-0.#2#3em #1}}
}
\def\rule_Buchstabe#1#2#3#4{%
 \suf_fix_def{#1}{%
  \scaled_rm_box{%
   \hbox{%
    #1%
    \hskip-0.#2em%
    \lower-0.#3ex\hbox{\vrule height1.#4ex width0.07em }%
   }%
   \hskip0.50em%
  }%
 }%
}
\I_Buchstabe B22
\rule_Buchstabe C51{34}
\I_Buchstabe D22
\I_Buchstabe E22
\I_Buchstabe F22
\rule_Buchstabe G{525}{081}4
\I_Buchstabe H22
\I_Buchstabe I20
\I_Buchstabe K22
\I_Buchstabe L20
\I_Buchstabe M{20em }{I\hskip-0.35}
\I_Buchstabe N{20em }{I\hskip-0.35}
\rule_Buchstabe O{525}{095}{45}
\I_Buchstabe P20
\rule_Buchstabe Q{525}{097}{47}
\I_Buchstabe R21 
\rule_Buchstabe U{45}{02}{54}
\suf_fix_def{Z}{\scaled_rm_box{Z\hskip-0.38em Z}}
\catcode`\"=12
\newcount\math_char_code
\def\suf_fix_math_chars_def#1{%
 \ifcat#1A
  \expandafter\math_char_code\expandafter=\suf_fix_fam
  \multiply\math_char_code by 256
  \advance\math_char_code by `#1
  \expandafter\mathchardef\csname#1\suf_fix\endcsname=\math_char_code
  \let\next=\suf_fix_math_chars_def
 \else
  \let\next=\relax
 \fi
 \next}
%
%
%
%
\def\font_fam_suf_fix#1#2 #3 {%
 \def\suf_fix{#2}
 \def\suf_fix_fam{#1}
 \suf_fix_math_chars_def #3.
}
\font_fam_suf_fix
 0rm
 ABCDEFGHIJKLMNOPQRSTUVWXYZabcdefghijklmnopqrstuvwxyz
\font_fam_suf_fix
 2scr
 ABCDEFGHIJKLMNOPQRSTUVWXYZ
\font_fam_suf_fix
 \slfam sl
 ABCDEFGHIJKLMNOPQRSTUVWXYZabcdefghijklmnopqrstuvwxyz
\font_fam_suf_fix
 \bffam bf
 ABCDEFGHIJKLMNOPQRSTUVWXYZabcdefghijklmnopqrstuvwxyz
\font_fam_suf_fix
 \ttfam tt
 ABCDEFGHIJKLMNOPQRSTUVWXYZabcdefghijklmnopqrstuvwxyz
\font_fam_suf_fix
 \ssfam
 ss
 ABCDEFGHIJKLMNOPQRSTUVWXYZabcdefgijklmnopqrstuwxyz
\catcode`\_=8
\def\Cdss{{\fam\ssfam
    \mkern 4.2 mu \mathchoice%
    {\vrule height 6.5pt depth -.55pt width 1pt}%
    {\vrule height 6.5pt depth -.57pt width 1pt}%
    {\vrule height 4.55pt depth -.28pt width .8pt}%
    {\vrule height 3.25pt depth -.19pt width .6pt}%
    \mkern -6.3mu C}}%
\def\Qdss{{\fam\ssfam
    \mkern 3.8 mu \mathchoice%
    {\vrule height 6.5pt depth -.67pt width 1pt}%
    {\vrule height 6.5pt depth -.7pt width 1pt}%
    {\vrule height 4.55pt depth -.44pt width .7pt}%
    {\vrule height 3.25pt depth -.3pt width .5pt}%
    \mkern -5.9mu Q}}%
%
%
%
%
%
%
%
%
%
%
%
%
%
%
\font\teneuf=eufm10 
\font\seveneuf=eufm7
\font\fiveeuf=eufm5
\newfam\euffam \def\euf{\fam\euffam\teneuf} 
\textfont\euffam=\teneuf \scriptfont\euffam=\seveneuf
\scriptscriptfont\euffam=\fiveeuf

       \def\dfr{{\euf d}}

       \def\gfr{{\euf g}}

       \def\xfr{{\euf x}}

\input xypic.tex

\def\Hom{{\rm Hom}}
\def\Ext{\mathop{\rm Ext}\nolimits}
\def\Rep{{\rm Rep}}
\def\Mod{{\rm Mod}}
\def\Vec{{\rm Vec}}
\def\Ind{{\rm Ind}}
\def\oiota{\widehat{\mathop{\otimes}\limits_{K,\iota}}}
\def\opi{\widehat{\mathop{\otimes}\limits_{K,\pi}}}

\parindent=0pt

\centerline{\bf Duality for admissible locally analytic representations}

\centerline{by Peter Schneider and Jeremy Teitelbaum}

\bigskip

{\bf 0. Introduction}

\medskip

In this paper we continue our development of the theory of locally
analytic representations of a $p$-adic locally analytic group $G$.  In
our earlier work [ST1,2] we constructed a certain abelian subcategory
of the category of modules over the locally analytic distribution
algebra $D(G,K)$, where $K$ is a complete discretely valued extension
of $\Qdss_p$. This subcategory of {\it coadmissible} modules is
contravariantly equivalent to the category of admissible locally
analytic representations by means of the functor sending a
representation $V$ to its strong dual $V'_b$. The smooth admissible
representations correspond to those coadmissible $D(G,K)$-modules
which are annihilated by the Lie algebra $\gfr$ of $G$.  Here, we
study the problem of constructing a contragredient functor on the
category of admissible locally analytic representations.  In fact, a
naive contragredient does not exist. As a best approximation, we
construct an involutive "duality" functor from the bounded derived
category of $D(G,K)$-modules with coadmissible cohomology to itself.
On the subcategory corresponding to complexes of smooth
representations, this functor induces the usual smooth contragredient
(with a degree shift). Although we construct our functor in general we
obtain its involutivity, for technical reasons, only in the case of
locally $\Qdss_p$-analytic groups.

The duality functor we construct in this paper is an extension to
coadmissible modules over $D(G,K)$ of the global duality $M\to
R\Hom_A(M,A)$ for an Auslander regular ring $A$.   Although $D(G,K)$
itself is not Auslander regular, it "almost" is.  More precisely it is
a free module over the distribution algebra $D(H,K)$ where $H$ is a
compact open subgroup of $G$.  In turn, $D(H,K)$ is, in the
terminology of [ST2], a Fr\'echet-Stein algebra, meaning that it is a
projective limit of a family of noetherian Banach algebras
$D_{r}(H,K)$ with flat transition maps. In Section 8 of [ST2] we
showed that, if the base field is $\Qdss_p$, these Banach algebras are
Auslander regular rings with global dimension bounded by the dimension
of $H$. Further, we showed that, for a coadmissible module $M$, the
modules $\Ext^{*}_{D(H,K)}(M,D(H,K))$ are coadmissible. The general
theory of coadmissible modules makes it possible to pass back and
forth between $D(H,K)$-modules and modules over the Auslander regular
Banach algebras $D_{r}(H,K)$. In addition, a kind of "Shapiro's Lemma"
makes it possible to pass back and forth between $D(G,K)$ and
$D(H,K)$.

The abelian category of coadmissible modules is filtered by a
generalized notion of grade or codimension.  Although the duality
functor on all coadmissible modules is expressed in terms of a derived
category, on the abelian subquotient categories corresponding to the
grade filtration, the duality functor is computed as a particular
$\Ext$-group.

We now briefly outline the paper. We begin by discussing the smooth
contragredient in the setting of coadmissible modules. Such
representations correspond to modules over the ring $D^{\infty}(G,K)$
of locally constant distributions, which is the quotient of $D(G,K)$
by the ideal generated by the Lie algebra $\gfr$. We then study the
dualizing modules $\Dscr^{\infty}_{K}(G)$ and $\Dscr_{K}(G)$ for
smooth and general locally analytic representations respectively.
These turn out to be the duals of the compactly supported smooth or
locally analytic functions on $G$.  Up to a twist which we suppress
here for simplicity, the duality functors on smooth and locally
analytic representations are given by $R\Hom_{A}(.,\Dscr)$ where $A$
is $D^{\infty}(G,K)$ or $D(G,K)$ and $\Dscr$ is
$\Dscr^{\infty}_{K}(G)$ or $\Dscr_{K}(G)$, though in the smooth case,
this is unnecessarily complicated since in fact
$\Dscr^{\infty}_{K}(G)$ is an injective module. In Section 3 we
compare the duality functors in the two categories. We use the
(trivial) deRham cohomology of the Lie group $G$ to compute (again,
suppressing some twists)
$R\Hom_{D(G,K)}(D^{\infty}(G,K),\Dscr_{K}(G))$ and from this
calculation we obtain the result that the duality on the locally
analytic representations induces the usual smooth duality, with a
degree shift, on the smooth representations. In Section 4, we
establish the involutive nature of the duality functor on the bounded
derived category of $D(G,K)$-modules with coadmissible cohomology.
Section 5 recalls how this duality respects the subquotient categories
of coadmissible modules of fixed codimension.  In the final section of
the paper we compute the duality functor for principal series
representations.  We show that a principal series representation
induced from a parabolic subgroup $P$ has a single nonvanishing
Ext-group in degree equal to the dimension of the subgroup $P$ which
is isomorphic to another principal series. It is interesting to
compare this with the similar result one obtains in the case of Verma
modules ([Kem]).

We thank Matthew Emerton for his suggestions regarding the
compatibility of the duality discussed in [ST2] and smooth duality.
The first author thanks the University of Illinois at Chicago and the
Clay Mathematics Institute for financial support. The second author
was supported by NSF grant DMS-0245410.

\medskip

{\bf 1. The problem}

\medskip

Let $\Qdss_p \subseteq L \subseteq K \subseteq \Cdss_p$ be complete
intermediate fields such that $L/\Qdss_p$ is finite and $K$ is
discretely valued. Throughout let $G$ be a locally $L$-analytic group
of dimension $d$ and let $\gfr$ denote its Lie algebra.

We let $D(G,K)$, resp. $D^\infty(G,K)$, be the $K$-algebra of
$K$-valued locally analytic, resp. locally constant, distributions on
$G$. The group $G$ embeds into $D(G,K)$ as well as into
$D^\infty(G,K)$ via the Dirac distributions $g \mapsto
\delta_g$. The universal enveloping algebra $U_K(\gfr) := U(\gfr)
\otimes_{L} K$ is naturally a subalgebra of $D(G,K)$. One has
$D^\infty(G,K) = D(G,K)/I_G(\gfr)$ where $I_G(\gfr)$ denotes the
closed two sided ideal of $D(G,K)$ generated by $\gfr$.

\smallskip

{\bf Remark 1.1:} {\it i. The group $G$ generates a dense subspace of
$D(G,K)$;

ii. the ideal $I_G(\gfr)$ is generated by $\gfr$ as a left ideal in
$D(G,K)$;

iii. $D(G,K) \otimes_{U(\gfr)} L = D^\infty(G,K)$;

iv. if $\gfr$ is semisimple then the ideal $I_G(\gfr)$ is idempotent.
}

Proof: i. [ST1] Lemma 3.1. ii. The identity $\delta_g
\xfr \delta_{g^{-1}} = {\rm ad}(g)(\xfr)$ for any $g \in G$ and $\xfr
\in \gfr$ together with i. imply that $I_G(\gfr)$ is the closure of the
left ideal generated by $\gfr$. On the other hand the left ideal
$D(G,K)\gfr$ certainly is finitely generated. If $G$ is compact it
therefore is closed by [ST2] Cor. 3.4.iv and Lemma 3.6. The general
case is reduced to the compact case by choosing a compact open
subgroup $H \subseteq G$ and observing the locally convex direct sum
decompositions
$$
D(G,K) = \bigoplus_{g \in G/H}
\delta_gD(H,K)\quad\hbox{and}\quad
D(G,K)\gfr = \bigoplus_{g \in G/H}
\delta_gD(H,K)\gfr\ .
$$
iii. This follows from ii. iv. This also follows from ii. since, by
[Dix] Remark 2.8.8, the ideal generated by $\gfr$ in $U(\gfr)$ is
idempotent.

\smallskip

We let $\Rep_K(G)$, resp. $\Rep_K^c(G)$ denote the category of locally
analytic $G$-representations on barrelled locally convex Hausdorff
$K$-vector spaces, resp. on $K$-vector spaces of compact type, with
continuous linear $G$-maps as morphisms (cf. [ST1] \S3). We recall:

-- On each $G$-representation in $\Rep_K(G)$ the $G$-action extends
uniquely to a separately continuous action of the algebra $D(G,K)$.

-- The category $\Rep_K^c(G)$ is closed with respect to the passage to
closed $G$-invariant subspaces and their corresponding quotients.

Dually we consider the abelian category $\Mscr_G$ of all (unital left)
$D(G,K)$-modules (in the algebraic sense), as well as the category
$\Mscr_G^{top}$ of all separately continuous $D(G,K)$-modules on
nuclear Fr\'echet spaces with continuous $D(G,K)$-modules maps as
morphisms. By functoriality, passing to the continuous dual is a
natural functor
$$\matrix{
\Rep_K(G) & \longrightarrow & \Mscr_G \cr
\hfill V & \longmapsto & V'\ . \hfill }
$$
It was proved in [ST1] Cor. 3.3 that the functor
$$
\matrix{
\Rep_K^c(G) & \mathop{\longrightarrow}\limits^{\simeq} & \Mscr_G^{top} \cr
 \hfill V & \longmapsto & V'_b \hfill}\leqno{(\ast)}
$$
of passing to the strong dual even is an anti-equivalence of
categories. (We silently use here, as at many other places in the
paper, the map $g \mapsto g^{-1}$ on $G$ to identify left and right
$D(G,K)$-modules.)

The objects in $\Rep_K^c(G)$ and in $\Mscr_G^{top}$ already as
topological vector spaces are of a quite different nature. This means
that the usual construction of a contragredient of a group
representation does not preserve the category $\Rep_K^c(G)$. The
question we want to address in this paper is whether there
nevertheless is a natural involutory functor on $\Rep_K^c(G)$ which
deserves to be considered as a replacement for the usual
contragredient. In fact we will restrict our attention to the full
subcategory $\Rep_K^a(G)$ of admissible $G$-representations in
$\Rep_K^c(G)$ as constructed in [ST2] \S6. We recall that in [ST2]
first a full abelian subcategory $\Cscr_G$ (which is entirely
algebraic in nature) of $\Mscr_G$, the category of coadmissible
$D(G,K)$-modules, is constructed. It then is shown (Lemma 6.1) that
$\Cscr_G$ embeds naturally and fully into the topological category
$\Mscr_G^{top}$. Finally $\Rep_K^a(G)$ is defined to be the preimage
of $\Cscr_G$ under the anti-equivalence $(\ast)$.

For the sake of completeness we also recall from [ST2] Prop. 6.4 that:

-- $\Rep_K^a(G)$ is an abelian category; kernel and image of a
morphism in $\Rep_K^a(G)$ are the algebraic kernel and image with the
subspace topology.

-- Any map in $\Rep_K^a(G)$ is strict and has closed image.

-- The category $\Rep_K^a(G)$ is closed with respect to the passage to
closed $G$-invariant subspaces.

An obvious condition which any solution of our problem should satisfy
comes about as follows. Let $\Rep_K^\infty(G)$, resp.
$\Rep_K^{\infty,a}(G)$ denote the category of smooth, resp. of
admissible-smooth, $G$-representations (over $K$) in the sense of
Jacquet-Langlands (cf. [Cas] \S2.1). Equipping  a $K$-vector space
with its finest locally convex topology induces a fully faithful
embedding
$$
\Rep_K^\infty(G) \longrightarrow \Rep_K(G)\ .
$$
It was proved in [ST2] Thm. 6.6 that this embedding restricts to a
fully faithful embedding
$$
\Rep_K^{\infty,a}(G) \longrightarrow \Rep_K^a(G)
$$
whose image consists precisely of those admissible locally analytic
$G$-representa- tions which are annihilated by the Lie algebra $\gfr$.
In the smooth theory one has the contragredient functor
$$
\matrix{
\Rep_K^\infty(G) & \longrightarrow & \Rep_K^\infty(G) \cr
  \hfill V & \longmapsto & \widetilde{V} \hfill }
$$
where the so called smooth dual $\widetilde{V}$ consists of all linear
forms in the full linear dual $V^\ast$ which are smooth, i.e., which
are fixed by some open subgroup of $G$. It restricts to an
anti-involution on the subcategory $\Rep_K^{\infty,a}(G)$ (cf. [Cas]
Prop. 2.1.10). Surely any construction which one might envisage on
$\Rep_K^a(G)$ should be compatible with this one.

It is easy to see that the $G$-action on any smooth $G$-representation
extends naturally to a $D^\infty(G,K)$-module structure. Because of
Remark 1.1.ii we may characterize $\Rep_K^{\infty,a}(G)$ as the full
subcategory of $\Rep_K^a(G)$ of those representations on which the
$D(G,K)$-module structure factorizes through a $D^\infty(G,K)$-module
structure. If $\Mscr_G^\infty$ denotes the abelian category of (unital
left) $D^\infty(G,K)$-modules then passing to the (full) linear dual
is a natural functor
$$
\matrix{
\Rep_K^\infty(G) & \longrightarrow & \Mscr_G^\infty \cr
  \hfill V & \longmapsto & V^\ast \hfill }
$$
which makes the diagram
$$
\xymatrix{
  \Rep_K^\infty(G) \ar[d]_{finest\ l.c.\ top.} \ar[rr]^{V \mapsto V^\ast} & &
  \Mscr_G^\infty \ar[d]^{\subseteq} \\
  \Rep_K(G) \ar[rr]^{V \mapsto V'} & & \Mscr_G   }
$$
commutative. We also may define the notion of a coadmissible
$D^\infty(G,K)$-module. This is based upon the following observation.

\smallskip

{\bf Remark 1.2:} {\it For any compact open subgroup $H \subseteq G$
we have:

i. $D^\infty(H,K)$ is a Fr\'echet-Stein algebra;

ii. $D^\infty(H,K)$ is coadmissible as a (left or right)
$D(H,K)$-module;

iii. a $D^\infty(H,K)$-module is coadmissible if and only if it is
coadmissible as a $D(H,K)$-module. }

Proof: i. Compare the proof of [ST2] Thm. 6.6.i. ii. This follows from
[ST2] Lemma 3.6. iii. This then is a special case of [ST2] Lemma 3.8.

\smallskip

We now define a $D^\infty(G,K)$-module to be coadmissible if it is so
as a $D^\infty(H,K)$-module for any compact open subgroup $H
\subseteq G$, and we let $\Cscr_G^\infty$ denote the full subcategory of
coadmissible $D^\infty(G,K)$-modules in $\Mscr_G^\infty$. We have
$$
\Cscr_G^\infty = \Cscr_G \cap \Mscr_G^\infty
$$
as well as the commutative diagram
$$
\xymatrix{
  \Rep_K^{\infty,a}(G) \ar[d] \ar[r]^{\quad\simeq} &
  \Cscr_G^\infty \ar[d] \\
  \Rep_K^a(G) \ar[r]^{\quad\simeq} & \Cscr_G\ .   }
$$
Finally we describe the smooth contragredient on the dual module side.
We introduce the $K$-vector space $C_c^\infty(G,K)$ of all $K$-valued
locally constant functions with compact support on $G$ and its linear
dual $\Dscr_K^\infty(G) := C_c^\infty(G,K)^\ast$. The left and right
translation actions of $G$ on $C_c^\infty(G,K)$ are smooth and
therefore extend to $D^\infty(G,K)$-module structures on
$C_c^\infty(G,K)$ and hence on $\Dscr_K^\infty(G)$. We fix a left
invariant Haar measure $\mu_{Haar} \in \Dscr_K^\infty(G)$. This allows
us to introduce, for any compact open subgroup $U \subseteq G$, the
function
$$
\epsilon_U := {\rm vol}_{\mu_{Haar}}(U)^{-1}\cdot {\rm char.\
function\ of}\ U
$$
in $C_c^\infty(G,K)$. One easily checks that
$$
\matrix{
C_c^\infty(G,K) & \longrightarrow & D^\infty(G,K) \hfill\cr
\hfill \varphi & \longmapsto & \mu_\varphi := \mu_{Haar}(\varphi\cdot
.) }
$$
is a left $G$-equivariant embedding. We also need the locally constant
modulus character $\delta_G : G \longrightarrow \Qdss^\times \subseteq
K^\times$ of $G$ given by
$$
\delta_G(g) = [gUg^{-1} : gUg^{-1} \cap U]/[U : gUg^{-1} \cap U]
$$
where $U$ is a fixed compact open subgroup of $G$. In the following we
will let $\delta_G$ also denote the one dimensional $K$-vector space
$K$ viewed as a $D^\infty(G,K)$-bimodule with $G$ acting trivially
from the left and via the character $\delta_G$ from the right. Finally
we need the functor
$$
\matrix{
\Mscr^\infty_G & \longrightarrow & \Rep^\infty_K(G) \hfill\cr
\hfill M & \longmapsto & M^{sm} := \{m \in M : Um = m\ \hbox{for some
open subgroup}\ U \subseteq G\}\ . }
$$

\smallskip

{\bf Lemma 1.3:} {\it For any module $M$ in $\Mscr^\infty_G$ we have
the natural $G$-equivariant isomorphism}
$$
\matrix{
M^{sm} & \mathop{\longrightarrow}\limits^{\cong} & (C_c^\infty(G,K)
\otimes_K \delta_G) \otimes_{D^\infty(G,K)} M \hfill\cr
\hfill m & \longmapsto & (\epsilon_U \otimes 1) \otimes m\qquad\hbox{if}\ Um
= m \ . \hfill }
$$

Proof: The inverse map is given by $(\varphi \otimes 1) \otimes m
\longmapsto \mu_{\varphi}m$.

\smallskip

{\bf Lemma 1.4:} {\it The diagram
$$
\xymatrix{
  \Rep_K^\infty(G) \ar[d]_{V \mapsto \widetilde{V}} \ar[rr]^{V \mapsto V^\ast} &
  & \Mscr_G^\infty \ar[d]^{\Hom_{D^\infty(G,K)}(.,\Dscr_K^\infty(G) \otimes_K \delta_G^\ast)} \\
  \Rep_K^\infty(G) \ar[rr]^{V \mapsto V^\ast} & & \Mscr_G^\infty   }
$$
is commutative. }

Proof: (Recall that $\Hom_{D^\infty(G,K)}(.,\Dscr_K^\infty(G)
\otimes_K \delta_G^\ast)$ is considered as a $D^\infty(G,K)$-module
via the right multiplication on the target.) We need to establish a
natural isomorphism between $\Hom_K((V^\ast)^{sm},K)$ and
$$
\matrix{
\Hom_{D^\infty(G,K)}(V^\ast,\Dscr_K^\infty(G) \otimes_K \delta_G^\ast)
\hfill\cr\cr \qquad\qquad
= \Hom_K((C^\infty_c(G,K) \otimes_K \delta_G)
\otimes_{D^\infty(G,K)} V^\ast,K)\ . }
$$
For this it certainly is sufficient to find a $G$-equivariant natural
isomorphism
$$
(V^\ast)^{sm} \cong (C^\infty_c(G,K) \otimes_K \delta_G)
\otimes_{D^\infty(G,K)} V^\ast
$$
before passing to the linear dual. But this is a special case of Lemma
1.3.

\smallskip

Because the smooth contragredient is an anti-involution on
$\Rep_K^{\infty,a}(G)$ it is clear from Lemma 1.4 that the functor
$\Hom_{D^\infty(G,K)}(.,\Dscr_K^\infty(G) \otimes_K \delta_G^\ast)$
restricts to an anti-involution on $\Cscr_G^\infty$.

\medskip

{\bf 2. Dualizing modules}

\medskip

In this section we will introduce an object analogous to
$\Dscr_K^\infty(G)$ in the locally analytic context. We define the
locally convex $K$-vector space $C_c^{an}(G,K)$ of $K$-valued locally
analytic functions with compact support on $G$ by picking a compact
open subgroup $H \subseteq G$ and by setting
$$
C_c^{an}(G,K) := \bigoplus_{g \in G/H} C^{an}(gH,K)
$$
with the locally convex direct sum topology on the right hand side. By
[Fea] 2.2.4 this definition is independent of the choice of $H$. In
fact, whenever $G = \mathop{\bigcup}\limits^\cdot_{i \in I} U_i$ is a
disjoint covering by compact open subsets $U_i$ we have
$$
C_c^{an}(G,K) = \bigoplus_{i \in I} C^{an}(U_i,K) \ .
$$
The space $C_c^{an}(G,K)$ is barrelled ([NFA] Ex. 3 after Cor. 6.16)
and the left and right translation actions of $G$ on it are locally
analytic and hence extend to separately continuous $D(G,K)$-module
structures. By functoriality we have corresponding $D(G,K)$-module
structures on the strong dual
$$
\Dscr_K(G) := C_c^{an}(G,K)'_b \ .
$$

\smallskip

{\it Remark 2.1:} i. In both actions each individual element in
$D(G,K)$ acts by a continuous endomorphism on $\Dscr_K(G)$.

ii. If $G$ is second countable then $G/H$ is countable for any compact
open subgroup $H \subseteq G$. Hence $C_c^{an}(G,K)$ is of compact
type by [ST1] Prop. 1.2(ii). In this situation [ST1] Cor. 3.3 says
that both $D(G,K)$-actions on $\Dscr_K(G)$ are separately continuous.

\smallskip

If $H \subseteq G$ is any compact open subgroup then, by [NFA] 9.10,
we have
$$
\Dscr_K(G) = \prod_{g \in G/H} D(gH,K)\ .
$$
The projection map $\ell_{G,H} : \Dscr_K(G) \longrightarrow D(H,K)$
onto the factor $D(H,K)$ in this decomposition is a canonical
$D(H,K)$-bimodule homomorphism. Moreover, we have
$$
m = (\delta_g \ell_{G,H}(\delta_{g^{-1}}m))_{g \in G/H}
\qquad\hbox{for any}\ m \in \Dscr_K(G). \leqno{(\ast)}
$$

\smallskip

{\bf Lemma 2.2:} {\it For any (left) $D(G,K)$-module $X$ the map
$$
\matrix{
\Hom_{D(G,K)}(X,\Dscr_K(G)) & \mathop{\longrightarrow}\limits^{\cong}
& \Hom_{D(H,K)}(X,D(H,K)) \cr
\hfill F & \longmapsto & \ell_{G,H} \circ F \hfill }
$$
is bijective and right $D(H,K)$-equivariant. }

Proof: By writing $X$ as the cokernel of a map between free
$D(G,K)$-modules it suffices to consider the case $X = D(G,K)$. Then
the left hand side is equal to $\Dscr_K(G)$ and the map becomes
$$
\matrix{
\Dscr_K(G) & \longrightarrow & \Hom_{D(H,K)}(D(G,K),D(H,K)) \cr
\hfill m & \longmapsto & \Lscr_m(\lambda) := \ell_{G,H}(\lambda m)\ .
\hfill}
$$
The injectivity of this map is immediate from $(\ast)$. For the
surjectivity we note that an element $\Lscr$ in the right hand side is
determined by its values $\Lscr(\delta_{g_0})$ for $g_0$ running over
a set of representatives for the cosets in $H\setminus G$. Define $m
:= (\delta_g\Lscr(\delta_{g^{-1}}))_{g \in G/H}$. Then
$$
\Lscr_m(\delta_{g_0}) = \ell_{G,H}(\delta_{g_0}m) =
\delta_{g_0}\delta_{g_0^{-1}}\Lscr(\delta_{g_0}) = \Lscr(\delta_{g_0})
$$
and hence $\Lscr_m = \Lscr$.

\smallskip

{\bf Proposition 2.3:} {\it For any bounded above complex $X^\cdot$ of
$D(G,K)$-modules and any compact open subgroup $H \subseteq G$ we have
a natural $D(H,K)$-equivariant isomorphism}
$$
\Ext^\ast_{D(G,K)}(X^\cdot,\Dscr_K(G)) \cong
\Ext^\ast_{D(H,K)}(X^\cdot,D(H,K))\ .
$$
Proof: Since $D(G,K)$ is free over $D(H,K)$ both sides can be computed
using a projective resolution of $X^\cdot$ as complex of
$D(G,K)$-modules. The assertion then is a consequence of the previous
lemma.

\smallskip

For the algebra $D^\infty(G,K)$ we have the following much stronger
fact.

\smallskip

{\bf Proposition 2.4:} {\it The left $D^\infty(G,K)$-module
$\Dscr^\infty_K(G)$ is injective.}

Proof: Fixing a compact open subgroup $H \subseteq G$ we obtain, in a
way completely analogous to the proof of Prop. 2.3, that
$$
\Ext^\ast_{D^\infty(G,K)}(Y,\Dscr^\infty_K(G)) \cong
\Ext^\ast_{D^\infty(H,K)}(Y,D^\infty(H,K))
$$
for any $D^\infty(G,K)$-module $Y$. This reduces us to proving that
$D^\infty(H,K)$ is a self-injective ring. But $D^\infty(H,K)$ is the
projective limit
$$
D^{\infty}(H,K) = \mathop{\lim\limits_{\longleftarrow}}\limits_{N}
K[H/N]
$$
of the algebraic group rings $K[H/N]$ where $N$ runs over the open
normal subgroups of $H$. It easily follows that $D^\infty(H,K)$ is a
direct product of finite dimensional simple $K$-algebras and as such
is self-injective by [Lam] Chap. I Cor. (3.11B).

\smallskip

The fact that such a vanishing result is not available over $D(G,K)$
forces us to work on the level of derived categories. As usual we let
$D^b(\Ascr)$, for any abelian category $\Ascr$, denote its bounded
derived category (which here is understood to be the derived category
of all complexes in $\Ascr$ with only finitely many nonzero cohomology
objects). Moreover, whenever $\Ascr_0
\subseteq \Ascr$ is a full abelian subcategory closed under extensions
we have the triangulated subcategory $D^b_{\Ascr_0}(\Ascr)$ of
$D^b(\Ascr)$ consisting of all those complexes whose cohomology
objects lie in $\Ascr_0$. For technical reasons we also will need the
``bounded below'' versions of these categories which, as usual, are
denoted by replacing the superscript ``b'' by ``+''

It is a simple consequence of [ST2] Remark 3.2 that the abelian
subcategories $\Cscr_G^\infty$ in $\Mscr_G^\infty$ and $\Cscr_G$ in
$\Mscr_G$ are closed under extensions. Hence we have the triangulated
subcategories $D^b_{\Cscr_G^\infty}(\Mscr_G^\infty)$ in
$D^b(\Mscr_G^\infty)$ and $D^b_{\Cscr_G}(\Mscr_G)$ in $D^b(\Mscr_G)$
available.

Since $\Dscr_K^\infty(G) \otimes_K \delta^\ast_G$ is a
$D^\infty(G,K)$-bimodule the functor
$$
R\Hom_{D^\infty(G,K)}(.,\Dscr_K^\infty(G) \otimes_K \delta^\ast_G) :
D^b(\Mscr_G^\infty)
\longrightarrow D^+(\Mscr_G^\infty)
$$
is well defined. Of course, by Prop. 2.4 (note that a twist preserves
injectivity), it is simply given by
$\Hom_{D^\infty(G,K)}(.,\Dscr_K^\infty(G) \otimes_K \delta^\ast_G)$.
The functors $.^\ast$ and $\widetilde{.}$ of passing to the full
linear and the smooth dual, respectively, are exact on smooth
representations and therefore pass directly to derived categories
where we denote them by the same symbols. As a consequence of Lemma
1.4 we then have the commutative diagram
$$
\xymatrix{
  D^b(\Rep_K^{\infty}(G)) \ar[d]_{\widetilde{.}} \ar[r]^(.6){.^\ast}
  & D^b(\Mscr_G^\infty) \ar[d]^{R\Hom_{D^\infty(G,K)}(.,\Dscr_K^\infty(G)
  \otimes_K \delta^\ast_G)} \\
  D^b(\Rep_K^{\infty}(G)) \ar[r]^(.6){.^\ast} & D^b(\Mscr_G^\infty)\ .
  }\leqno{(2.5)}
$$
Correspondingly, for a certain twist $\Dscr_K(G) \otimes_K \dfr_G$ of
the $D(G,K)$-bimodule $\Dscr_K(G)$, to be defined in the next section,
we have the functor
$$
R\Hom_{D(G,K)}(.,\Dscr_K(G) \otimes_K \dfr_G) : D^b(\Mscr_G)
\longrightarrow D^+(\Mscr_G)
$$
In the next section we will establish the relation between these two
$R\Hom$-functors.

\medskip

{\bf 3. Lie algebra cohomology}

\medskip

We first recall the standard complexes which compute Lie algebra
(co)homology for an arbitrary $\gfr$-module $X$ (over $K$) (cf. [CE]
Chap. XIII\S8). Let $\bigwedge^\cdot\gfr$ denote, as usual, the
exterior algebra over $\gfr$. The Lie algebra cohomology
$H^\ast(\gfr,X)$ is the cohomology of the (bounded) complex which in
degree $q$ is
$$
C^q(\gfr;X) := \Hom_L(\bigwedge^q \gfr,X)
$$
and whose differential is given by
$$
\eqalign{
dc(\xfr_0,\ldots,\xfr_q) := & \sum_{s < t}
(-1)^{s+t+1}c([\xfr_s,\xfr_t],\xfr_1,\ldots,\widehat{\xfr_s},\ldots,
\widehat{\xfr_t},\ldots,\xfr_q) \cr
& + \sum_s (-1)^{s+1}\xfr_s
c(\xfr_1,\ldots,\widehat{\xfr_s},\ldots,\xfr_q)\ . \cr }
$$
Correspondingly the Lie algebra homology $H_\ast(\gfr,X)$ is the
homology of the complex which in degree $q$ is
$$
C_q(\gfr;X) := \bigwedge^q \gfr \otimes_L X
$$
and whose differential is given by
$$
\eqalign{
\partial(\xfr_1 \wedge\ldots\wedge \xfr_q \otimes x) := & \sum_{s < t}
(-1)^{s+t}[\xfr_s,\xfr_t] \wedge
\xfr_1 \wedge\ldots\wedge \widehat{\xfr_s} \wedge\ldots\wedge
\widehat{\xfr_t} \wedge\ldots\wedge \xfr_q \otimes x \cr
& + \sum_s (-1)^{s+1} \xfr_1
\wedge\ldots\wedge \widehat{\xfr_s} \wedge\ldots\wedge \xfr_q
\otimes \xfr_sx \ . \cr }
$$

The starting point of our investigation is the following basic fact.

\smallskip

{\bf Proposition 3.1:} {\it With respect to the natural either left or
right $\gfr$-module structure on $D(G,K)$ we have
$$
H_q(\gfr,D(G,K)) = \left\{\matrix{D^\infty(G,K) & \hbox{if}\ q = 0,
\cr 0 \hfill & \hbox{if}\ q > 0\cr}\right.
$$
and
$$
H_q(\gfr,\Dscr_K(G)) = \left\{\matrix{\Dscr^\infty_K(G) & \hbox{if}\ q
= 0,
\cr 0 \hfill & \hbox{if}\ q > 0.\cr}\right.
$$
}

Proof: By symmetry it suffices to consider the left $\gfr$-module
structure. Furthermore, for any compact open subgroup $H
\subseteq G$ we have the $\gfr$-invariant decompositions $D(G,K) =
\bigoplus_{g \in H\setminus G} D(H,K)\delta_g$ and $\Dscr_K(G) =
\prod_{g \in H\setminus G} D(H,K)\delta_g$. Since Lie algebra homology
commutes with arbitrary direct sums and direct products we may in fact
assume that $G$ is compact. We consider now the deRham complex
$$
0 \longrightarrow C^{an}(G,K) = A^0(G,K)
\mathop{\longrightarrow}\limits^d A^1(G,K)
\mathop{\longrightarrow}\limits^d \ldots
\mathop{\longrightarrow}\limits^d A^d(G,K) \longrightarrow 0
$$
of $K$-valued global locally analytic differential forms on the
locally analytic manifold $G$. By the usual Poincar\'e lemma it is an
exact resolution of the space $C^\infty(G,K)$ of locally constant
functions on $G$. Since the tangent bundle $TG = \gfr \times G$ on $G$
is trivial we have $A^q(G,K) = \Hom_L(\bigwedge^q\gfr,C^{an}(G,K))$.
This identifies (up to a sign) the deRham complex with the
cohomological standard complex for the $\gfr$-module $C^{an}(G,K)$
(cf. [BW] VII.1.1) and proves that $H^0(\gfr,C^{an}(G,K)) =
C^\infty(G,K)$ and $H^q(\gfr,C^{an}(G,K)) = 0$ for $q > 0$. But
because of the reflexivity of the vector space of compact type
$C^{an}(G,K)$ we can go one step further and have
$$
A^q(G,K) = \Hom_L(\bigwedge^q\gfr,C^{an}(G,K)) =
\Hom_K^{cont}(\bigwedge^q\gfr \otimes_L D(G,K),K)
$$
where $\Hom_K^{cont}$ on the right hand side refers to the continuous
linear forms. Hence the deRham complex in fact is the continuous dual
of the homological standard complex for the $\gfr$-module $D(G,K)$.
This homological standard complex is a complex of finitely generated
free right $D(G,K)$-modules. Its differentials therefore are
continuous and strict maps ([ST2] paragraph between Lemma 3.6 and
Prop. 3.7). Using the Hahn-Banach theorem we now see that the
continuous dual of $H_\ast(\gfr,D(G,K))$ is the cohomology
$H^\ast(\gfr,C^{an}(G,K))$ of the deRham complex which we computed
already.

\smallskip

By Remark 1.1.iii we have
$$
D(G,K) \otimes_{U_K(\gfr)} K = D^\infty(G,K)\ .
$$
It follows that the diagram of functors
$$
\xymatrix{
  \Mscr_G \ar[d]_{\rm forget} \ar[rrrr]^{\Hom_{D(G,K)}(D^\infty(G,K),.)}
                    & & & & \Mscr^\infty_G \ar[d]^{\rm forget} \\
  \Mod(U_K(\gfr)) \ar[rrrr]^{\Hom_{U_K(\gfr)}(K,.)} & & & & \Vec_K   }
$$
is commutative where $\Mod(U_K(\gfr))$ denotes the category of (unital
left) $U_K(\gfr)$-modules and $\Vec_K$ the category of $K$-vector
spaces . If we choose a projective resolution $P.$ of $K$ in
$\Mod(U_K(\gfr))$ then Prop. 3.1 implies that $D(G,K)
\otimes_{U_K(\gfr)} P.$ is a projective resolution of $D^\infty(G,K)$
in $\Mscr_G$. Choosing now also an injective resolution $I^\cdot$ of
an $X$ in $\Mscr_G$ we obtain the following sequence of identities
$$
\matrix{
R\Hom_{U_K(\gfr)}(K,{\rm forget}(X)) & \sim & \Hom_{U_K(\gfr)}(P.,{\rm
forget}(X))\hfill\cr
 & \sim & \Hom_{D(G,K)}(D(G,K) \otimes_{U_K(\gfr)}P.,X)\hfill\cr
 & \sim & {\rm forget}(\Hom_{D(G,K)}(D^\infty(G,K),I^\cdot))\hfill\cr
 & \sim & {\rm forget}(R\Hom_{D(G,K)}(D^\infty(G,K),X))\hfill }
$$
in the derived category of $\Vec_K$. Hence we have
$$
\Ext^q_{D(G,K)}(D^\infty(G,K),X) = \Ext^q_{U_K(\gfr)}(K,X) =
H^q(\gfr,X)
$$
for any $q \geq 0$ (we suppress the forgetful functors in the
notation) and, in particular, $\Ext^q_{D(G,K)}(D^\infty(G,K),X) =
H^q(\gfr,X) = 0$ for $q > d$. Therefore the corresponding total right
derived functors form a commutative diagram on the level of bounded
derived categories:
$$
\xymatrix{
 D^b(\Mscr_G) \ar[d]_{\rm forget} \ar[rrrr]^{R\Hom_{D(G,K)}(D^\infty(G,K),.)}
                    & & & & D^b(\Mscr^\infty_G) \ar[d]^{\rm forget} \\
  D^b(\Mod(U_K(\gfr))) \ar[rrrr]^{R\Hom_{U_K(\gfr)}(K,.)} & & & & D^b(\Vec_K)   }
$$
On the other hand, since the functor $\Hom_{D(G,K)}(D^\infty(G,K),.) :
\Mscr_G \longrightarrow \Mscr^\infty_G$ preserves injective objects,
we obviously have the adjointness relations
$$
\Hom_{ D^b(\Mscr_G)}(Y^\cdot,X^\cdot) =
\Hom_{D^b(\Mscr^\infty_G)}(Y^\cdot,R\Hom_{D(G,K)}(D^\infty(G,K),X^\cdot))
$$
and
$$
\matrix{
R\Hom_{D(G,K)}(Y^\cdot,X^\cdot) = \hfill \cr \qquad
R\Hom_{D^\infty(G,K)}(Y^\cdot,R\Hom_{D(G,K)}(D^\infty(G,K),X^\cdot))}
\leqno{(+)}
$$
in $D^b(\Vec_K)$ and in particular
$$
\Ext^\ast_{D(G,K)}(Y^\cdot,X^\cdot) =
\Ext^\ast_{D^\infty(G,K)}(Y^\cdot,R\Hom_{D(G,K)}(D^\infty(G,K),X^\cdot))
$$
for any $Y^\cdot$ in $D^b(\Mscr^\infty_G)$ and any $X^\cdot$ in
$D^b(\Mscr_G)$.

Unfortunately it does not seem to be the case that
$R\Hom_{D(G,K)}(D^\infty(G,K),.)$ ``restricts'' to a functor from
$D^b(\Cscr_G)$ into $D^b(\Cscr_G^\infty)$.

But there is a little bit more which can be said about this functor. A
particular finitely generated free resolution of $K$ in
$\Mod(U_K(\gfr))$ is given by the homological standard complex
$$
U_K(\gfr) \otimes_L \bigwedge^d \gfr
\mathop{\longrightarrow}\limits^{\partial} \ldots
\mathop{\longrightarrow}\limits^{\partial}
U_K(\gfr) \otimes_L \bigwedge^1 \gfr
\mathop{\longrightarrow}\limits^{\partial}
U_K(\gfr) \otimes_L \bigwedge^0 \gfr \longrightarrow K
$$
for $U_K(\gfr)$ as a right $\gfr$-module. Base extending this complex
to $D(G,K)$, by Prop. 3.1, is exact and therefore provides the
finitely generated free ``standard'' resolution
$$
D(G,K) \otimes_L \bigwedge^d \gfr
\mathop{\longrightarrow}\limits^{\partial} \ldots
\mathop{\longrightarrow}\limits^{\partial}
D(G,K) \otimes_L \bigwedge^0 \gfr \longrightarrow D^\infty(G,K)
\leqno{(\ast)}
$$
of $D^\infty(G,K)$ in $\Mscr_G$. We recall that the differential
$\partial$ is given by
$$
\eqalign{
\partial(\lambda \otimes \xfr_1 \wedge\ldots\wedge \xfr_q) = & \sum_{s < t}
(-1)^{s+t} \lambda \otimes [\xfr_s,\xfr_t] \wedge
\xfr_1 \wedge\ldots\wedge \widehat{\xfr_s} \wedge\ldots\wedge
\widehat{\xfr_t} \wedge\ldots\wedge \xfr_q \cr
& + \sum_s (-1)^{s+1} \lambda\xfr_s \otimes \xfr_1
\wedge\ldots\wedge \widehat{\xfr_s} \wedge\ldots\wedge \xfr_q
\ . \cr }
$$
We claim that the resolution $(\ast)$ carries a natural commuting
right $D(G,K)$-module structure which on $D^\infty(G,K)$ is the
obvious one. The group $G \subseteq D(G,K)$ obviously acts on $D(G,K)
\otimes_L \bigwedge^q \gfr$ from the right by
$$
(\lambda \otimes \xfr_1 \wedge\ldots\wedge \xfr_q)\delta_g :=
\lambda\delta_g \otimes {\rm ad}(g^{-1})(\xfr_1) \wedge\ldots\wedge
{\rm ad}(g^{-1})(\xfr_q) \ .
$$
One easily checks that the differentials $\partial$ respect this
$G$-action.

\smallskip

{\bf Lemma 3.2:} {\it The above $G$-action on $(\ast)$ extends
uniquely to a separately continuous $D(G,K)$-module structure.}

Proof: See the Appendix. The equivariance of the differentials follows
by continuity.

\smallskip

The usual argument with double complexes now shows that, for any
bounded complex $X^\cdot$ of $D(G,K)$-modules, we indeed have
$$
R\Hom_{D(G,K)}(D^\infty(G,K),X^\cdot) \sim \Hom^\cdot_{D(G,K)}(D(G,K)
\otimes \bigwedge^\cdot \gfr, X^\cdot)
$$
in $D^b(\Mscr_G)$. This will enable us to reinterprete this functor in
a completely different way. We write in the following $\Delta_G :=
\bigwedge^d \gfr \otimes_L K$ considered as a one dimensional locally
analytic $G$-representation. In particular, $\Delta_G$ is a one
dimensional $D(G,K)$-module given by an algebra homomorphism which, by
abuse of notation, we also write $\Delta_G : D(G,K) \longrightarrow
K$. In the most interesting case where $G$ is open in the group of
$K$-valued points of a connected reductive $K$-group this homomorphism
$\Delta_G$ in fact is trivial. But there is no point in restricting to
this case in the following. As will be shown in the Appendix, the
tensor product $\Delta_G
\otimes_K X$ with any other $D(G,K)$-module $X$ is defined as a
$D(G,K)$-module. Since $\Delta_G$ is one dimensional a more explicit
way to describe the tensor product $\Delta_G \otimes_K X$ is to say
that it is the pull back of the $D(G,K)$-module $X$ via the algebra
automorphism
$$
\matrix{
\alpha_{\Delta_G} : D(G,K) & \longrightarrow & D(G,K) \hfill\cr
            \hfill \lambda & \longmapsto &
            \alpha_{\Delta_G}(\lambda)(\psi) :=
            \lambda(\Delta_G\psi)\ . }
$$
By our above standard resolution the $D^\infty(G,K)$-module
$\Ext^d_{D(G,K)}(D^\infty(G,K),\break
\Delta_G \otimes_K X)$ is the cokernel of the map
$$
\xymatrix{
\Hom_{D(G,K)}(D(G,K)
\otimes_L \bigwedge^{d-1} \gfr, \Delta_G \otimes_K X)
\ar[d]^{\partial^\ast} \\
\Hom_{D(G,K)}(D(G,K) \otimes_L \bigwedge^d \gfr,
\Delta_G \otimes_K X)\ . }
$$
For the target we have the natural isomorphism
$$
\matrix{
X & \mathop{\longrightarrow}\limits^{\cong} & \Hom_{D(G,K)}(D(G,K)
\otimes_L \bigwedge^d \gfr, \Delta_G \otimes_K X) \cr
x & \longmapsto & [\lambda \otimes \delta \mapsto \lambda(\delta
\otimes x)] \ .\hfill }
$$
It in fact is an isomorphism of $D(G,K)$-modules where the left
$D(G,K)$-module structure on the right hand side is induced by the
right $D(G,K)$-module structure on $D(G,K)
\otimes_L \bigwedge^d \gfr$; this is easily checked by using
the above explicit description of the tensor product modules. To
compute the cokernel we fix an $L$-basis $\xfr_1,\ldots,\xfr_d$ of
$\gfr$, and observe right away the elementary identity
$$
\matrix{
\sum_{i < j} (-1)^{i+j} [\xfr_i,\xfr_j] \wedge
\xfr_1 \wedge\ldots\wedge \widehat{\xfr_i} \wedge\ldots\wedge
\widehat{\xfr_j} \wedge\ldots\wedge \xfr_d \hfill\cr\cr
\qquad\qquad = \sum_s (-1)^{s} {\rm tr}({\rm Ad}(\xfr_s)) \xfr_1
\wedge\ldots\wedge \widehat{\xfr_s} \wedge\ldots\wedge \xfr_d
\ . }
$$
If write an element $f \in \Hom_{D(G,K)}(D(G,K)
\otimes_L \bigwedge^{d-1} \gfr, \Delta_G \otimes_K X) =\break
\Hom_L(\bigwedge^{d-1} \gfr ,\Delta_G \otimes_K X)$ as $f(.) =
\xfr_1\wedge\ldots\wedge\xfr_d \otimes f_X(.)$ then it is
straightforward from the above identity that
$$
\partial^\ast f(\xfr_1\wedge\ldots\wedge\xfr_d) =
\xfr_1\wedge\ldots\wedge\xfr_d \otimes \sum_s (-1)^{s+1}\xfr_s f_X(\xfr_1
\wedge\ldots\wedge \widehat{\xfr_s} \wedge\ldots\wedge \xfr_d) \ .
$$
In other words, under the above identification of the target with $X$
the image $\partial^\ast f$ becomes the element $\sum_s
(-1)^{s+1}\xfr_s f_X(\xfr_1
\wedge\ldots\wedge \widehat{\xfr_s} \wedge\ldots\wedge \xfr_d)$ in
$X$. Since the values of $f_X$ can be chosen arbitrarily this means
that the cokernel identifies naturally with $X/\gfr X = X/D(G,K)\gfr X
= D^\infty(G,K) \otimes_{D(G,K)} X$. Hence we have established a
natural isomorphism of $D^\infty(G,K)$-modules
$$
\Ext^d_{D(G,K)}(D^\infty(G,K),
\Delta_G \otimes_K X) \cong D^\infty(G,K) \otimes_{D(G,K)} X \ .
\leqno{(1)}
$$
We repeat that
$$
\Ext^q_{D(G,K)}(D^\infty(G,K), \Delta_G \otimes_K X) = 0
\qquad\hbox{for}\ q > d\ . \leqno{(2)}
$$
Let $V$ be any $K$-vector space and consider $D(G,K)
\otimes_K V$ as a left $D(G,K)$-module in the obvious way. Since
$$
\matrix{
{\rm forget}(\Ext^\ast_{D(G,K)}(D^\infty(G,K), \Delta_G \otimes_K
D(G,K) \otimes_K V)) \cr\cr
 \qquad\qquad = H^\ast(\gfr,\Delta_G \otimes_K D(G,K) \otimes_K V)
 \hfill\cr\cr
 \qquad\qquad = H_{d-\ast}(\gfr,D(G,K) \otimes_K V) =
 H_{d-\ast}(\gfr,D(G,K)) \otimes_K V \hfill\cr\cr
 \qquad\qquad = 0 \hfill }
$$
for $\ast \neq d$, where the second, resp. third, identity comes from
[CE] XIII Ex. 15, resp. from Prop. 3.1, we have
$$
\Ext^q_{D(G,K)}(D^\infty(G,K), \Delta_G \otimes_K D(G,K) \otimes_K V) = 0
\qquad\hbox{for}\ q \neq d\ .\leqno{(3)}
$$
Finally we observe that any $D(G,K)$-module has a resolution by free
$D(G,K)$-modules of the form $D(G,K) \otimes_K V$ (use inductively the
natural surjection $D(G,K) \otimes_K X \rightarrow X$). In this
situation, with (1) -- (3), general homological algebra ([Har] I.7.4)
tells us that $\Ext^\ast_{D(G,K)}(D^\infty(G,K),.)$ is the left
derived functor of $\Ext^d_{D(G,K)}(D^\infty(G,K),\Delta_G \otimes_K
.) \cong D^\infty(G,K) \otimes_{D(G,K)} .$ This establishes the
following fact.

\smallskip

{\bf Proposition 3.3:} {\it There is a natural isomorphism of
$D^\infty(G,K)$-modules
$$
\Ext^\ast_{D(G,K)}(D^\infty(G,K),\Delta_G \otimes_K X)
\cong {\rm Tor}^{D(G,K)}_{d-\ast} (D^\infty(G,K),X)
$$
for any $D(G,K)$-module $X$.}

\smallskip

For our purpose of understanding duality the following special case of
these considerations is the most interesting. We introduce the one
dimensional tensor product
$$
\dfr_G := \Delta_G \otimes_K \delta^\ast_G
$$
viewed as a $D(G,K)$-bimodule with $D(G,K)$ acting trivially from the
right and through the product character $\Delta_G \cdot \delta_G$ from
the left.

\smallskip

{\it Remark 3.4:} By [B-GAL] Chap. III\S3.16 Cor. to Prop. 55 we have
$\delta_G = |\Delta_G|_L$ where $|.|_L$ denotes the normalized
absolute value of the field $L$.

\smallskip

{\bf Proposition 3.5:} {\it $R\Hom_{D(G,K)}(D^\infty(G,K),\Dscr_K(G)
\otimes_K \dfr_G)$ is naturally quasi-isomorphic to
$\Dscr^\infty_K(G) \otimes_K \delta^\ast_G$ concentrated in degree
$d$.}

Proof: Since $\gfr$ acts trivially on $\delta^\ast_G$ a computation
analogous to the one between formulae $(2)$ and $(3)$ and based upon
the second part of Prop. 3.1 shows that
$$
\Ext^q_{D(G,K)}(D^\infty(G,K), \Delta_G \otimes_K \Dscr_K(G) \otimes_K \delta^\ast_G) = 0
\qquad\hbox{for}\ q \neq d\ ,
$$
i.e., that the complex in question is concentrated in degree $d$. Its
cohomology in degree $d$, by formula $(1)$, is naturally isomorphic to
$D^\infty(G,K) \otimes_{D(G,K)} (\Dscr_K(G) \otimes_K \delta^\ast_G)$.
Using Remark 1.1 we obtain
$$
\matrix{
D^\infty(G,K) \otimes_{D(G,K)} (\Dscr_K(G) \otimes_K \delta^\ast_G)
\hfill\cr\cr \qquad\qquad
= (\Dscr_K(G) \otimes_K \delta^\ast_G)/\gfr(\Dscr_K(G) \otimes_K
\delta^\ast_G) \hfill\cr\cr \qquad\qquad
= (\Dscr_K(G) \otimes_K \delta^\ast_G)/(\gfr\Dscr_K(G) \otimes_K
\delta^\ast_G) \hfill\cr\cr \qquad\qquad
= (\Dscr_K(G)/\gfr\Dscr_K(G)) \otimes_K \delta^\ast_G \hfill\cr\cr
\qquad\qquad
= \Dscr^\infty_K(G) \otimes_K \delta^\ast_G \hfill}
$$
with the last identity coming from Prop. 3.1.

\smallskip

{\bf Corollary 3.6:} {\it The diagram
$$
\xymatrix{
  D^b(\Mscr^\infty_G) \ar[d]_{R\Hom_{D^\infty(G,K)}(.,\Dscr^\infty_K(G) \otimes_K \delta^\ast_G[-d])}
      \ar[r]^{\rm can} &  D^b(\Mscr_G) \ar[d]^{R\Hom_{D(G,K)}(.,\Dscr_K(G) \otimes_K \dfr_G)} \\
  D^b(\Mscr^\infty_G) \ar[r]^{\rm can} & D^+(\Mscr_G) }
$$
is commutative.}

Proof: This follows from Prop. 3.5 and the special case
$$
\matrix{
R\Hom_{D(G,K)}(Y^\cdot,\Dscr_K(G) \otimes_K \dfr_G) = \hfill \cr
\qquad R\Hom_{D^\infty(G,K)}(Y^\cdot,R\Hom_{D(G,K)}(D^\infty(G,K),
\Dscr_K(G) \otimes_K \dfr_G))}
$$
of our earlier adjointness relation $(+)$ provided we show that the
latter holds true actually in $D^+(\Mscr_G)$. For this we consider the
$D(G,K)-D^\infty(G,K)$-bimodule $D^\infty(G,K)$ as a module over the
ring $D(G,K) \otimes_K D^\infty(G,K)$ and fix a resolution $F.
\mathop{\longrightarrow}\limits^{\sim} D^\infty(G,K)$ by free
$D(G,K) \otimes_K D^\infty(G,K)$-modules. The point is that this
remains a free resolution if only considered as (left)
$D(G,K)$-modules or (right) $D^\infty(G,K)$-modules, respectively. If
in addition we choose a projective resolution $P^\cdot
\mathop{\longrightarrow}\limits^{\sim} Y^\cdot$ in $\Mscr^\infty_G$
then the right hand side of the above identity is computed by the
complex of $D(G,K)$-modules
$$
\Hom^\cdot_{D^\infty(G,K)}(P^\cdot,\Hom_{D(G,K)}(F.,
\Dscr_K(G) \otimes_K \dfr_G)) \ .
$$
But this complex is equal to the complex
$$
\Hom_{D(G,K)}(F. \otimes_{D^\infty(G,K)} P^\cdot,
\Dscr_K(G) \otimes_K \dfr_G) \ .
$$
Since $F. \otimes_{D^\infty(G,K)} P^\cdot
\mathop{\longrightarrow}\limits^{\sim} D^\infty(G,K)
\otimes_{D^\infty(G,K)} Y^\cdot = Y^\cdot$ is a free resolution of
$Y^\cdot$ in $\Mscr_G$ the second Hom-complex computes the left hand
side of our identity.

\smallskip

{\bf Corollary 3.7:} {\it The diagram
$$
\xymatrix{
  D^b(\Rep_K^{\infty}(G)) \ar[d]_{\widetilde{.}\,[-d]} \ar[r]^{\ .^\ast}
  & D^b(\Mscr_G) \ar[d]^{R\Hom_{D(G,K)}(.,\Dscr_K(G)
  \otimes_K \dfr_G)} \\
  D^b(\Rep_K^{\infty}(G)) \ar[r]^{\ .^\ast} & D^+(\Mscr_G)\ .
  }
$$
is commutative.}

Proof: Combine (2.5) and Cor. 3.6.

\medskip

{\bf Appendix:} The ``comultiplication'' on $D(G,K)$

\medskip

In fact, there is no true comultiplication on $D(G,K)$ as will become
clear presently. In particular, for two (left) $D(G,K)$-modules $X_1$
and $X_2$ the tensor product $X_1 \otimes_K X_2$ does not carry in
general a natural $D(G,K)$-module structure. But we do have the
following. Let $L \subseteq K$ be two extension fields of $\Qdss_p$
such that $L/\Qdss_p$ is finite and $K$ is spherically complete, and
let $M$ and $N$ be two paracompact locally $L$-analytic manifolds.
Furthermore let, as usual, the symbols $(\oiota)\
\mathop{\otimes}\limits_{K,\iota}$ and $(\opi)\
\mathop{\otimes}\limits_{K,\pi}$ stand for the (completed) inductive
and projective, respectively, tensor product of two locally convex
$K$-vector spaces (cf. [NFA] \S17).

\smallskip

{\bf Lemma A.1:} {\it If $N$ is compact then we have}
$$
C^{an}(M \times N,K) = C^{an}(M,C^{an}(N,K))\ .
$$
Proof: This is an exercise, which we leave to the reader, in the
construction of the locally convex vector spaces $C^{an}(.,.)$ based
upon the observation (cf. [Fea] 2.1.2) of a corresponding fact for
Banach space valued power series.

\smallskip

{\bf Proposition A.2:} {\it Suppose that $M$ is compact, and that $V$
is a $K$-vector space of compact type; we then have}
$$
C^{an}(M,K)\, \opi\, V = C^{an}(M,V)\ .
$$
Proof: Apply [Fea] 2.4.3 to the bilinear map $C^{an}(M,K) \times V \to
C^{an}(M,V)$ to see that the map is continuous and has dense image. If
we show that $C^{an}(M,V)$ induces the projective tensor topology on
$C^{an}(M,K) \otimes_K V$, then by completeness we have a topological
isomorphism.  By the integration theorem ([ST1] Thm. 2.2), we have a
linear isomorphism $C^{an}(M,V) \to {\cal L}(D(M,K),V)$.  At the same
time, using the fact that $V$ is the strong dual of a Fr\'{e}chet
space, we see by the discussion before Prop. 20.13 of [NFA] that we
have
$$
C^{an}(M,K)\, \opi\, V
\mathop{\longrightarrow}\limits^{\cong}
D(M,K)' \, \opi\, V
\mathop{\longrightarrow}\limits^{\cong}
{\cal L}_{b}(D(M,K),V)\ .
$$
Therefore the map $C^{an}(M,K)\, \opi\, V \to C^{an}(M,V)$ is
bijective, as well as continuous.   To complete the proof, we wish to
apply the open mapping theorem and for this ([NFA] 8.8) it suffices to
show that $C^{an}(M,K)\, \opi\, V$ is a direct limit of Banach spaces;
we in fact show that it is of compact type. Indeed, $D(M,K)$ and
$V'_{b}$ being reflexive Fr\'{e}chet spaces, are nuclear ([NFA] 20.7)
and therefore ([NFA] 19.11 and 20.4, 20.13, 20.14) $D(M,K)\,
\opi\, V'_{b}$ is again nuclear and reflexive, and
$$
(D(M,K)\, \opi\, V'_{b})'_{b} = D(M,K)'_b\,
\opi\, V = C^{an}(M,K)\, \opi\, V \ .
$$
As the strong dual of a nuclear Frechet space, $C^{an}(M,K)\,
\opi\, V$ is of compact type by [ST1] Thm. 1.3. (For a different proof
compare [Eme] Prop. 2.1.28.)

\smallskip

{\bf Proposition A.3:} $D(M \times N,K) = D(M,K)\, \oiota\, D(N,K)$.

Proof: We first assume that $M$ and $N$ are compact. Recall that the
Fr\'echet space $D(M,K)$ is the strong dual of the reflexive space of
compact type $C^{an}(M, K)$ (cf. [ST1] \S2). Combining A.1 and A.2 we
have
$$
C^{an}(M \times N,K) = C^{an}(M,K)\, \opi\, C^{an}(N,K)\ .
$$
Using reflexivity in addition we therefore obtain
$$
\eqalign{
D(M \times N,K) & = C^{an}(M \times N,K)'_b \cr
                    & = [C^{an}(M,K)\, \opi\, C^{an}(N,K)]'_b \cr
                    & = [D(M,K)'_b\, \opi\, D(N,K)'_b]'_b \ . \cr
                    }
$$
By [NFA] 20.13 and 20.14 the last term is equal to
$$
([D(M,K)\, \opi\, D(N,K)]'_b)'_b = D(M,K)\, \opi\, D(N,K)\ .
$$
Since for Fr\'echet spaces inductive and projective tensor product
topology coincide ([NFA] 17.6) our assertion
$$
D(M \times N,K) = D(M,K)\, \oiota\, D(N,K)
$$
follows in the compact case.

In the general case the assumption on paracompactness implies ([Sch]
Kap. II Satz 8.6) that there are decompositions $M =
\mathop{\bigcup}\limits^\cdot_{i \in I} M_i$ and $N =
\mathop{\bigcup}\limits^\cdot_{j \in J} N_j$ into disjoint unions of
open compact submanifolds $M_i$ and $N_j$, respectively. Since the
completed inductive tensor product commutes with arbitrary locally
convex direct sums by an obvious nonarchimedean version of [Gro] Chap.
I\S3.1 Prop. 14.I we may deduce the general case now as follows:
$$
\eqalign{
D(M \times N,K) & = \bigoplus_{i \in I} \bigoplus_{j \in J}
                             D(M_i \times N_j,K) \cr
                     & = \bigoplus_{i,j}
                             [D(M_i,K)\, \oiota\, D(N_j,K)] \cr
                     & = [\bigoplus_{i}  D(M_i,K)]\, \oiota\,
                             [\bigoplus_{j} D(N_j,K)] \cr
                     & = D(M,K)\, \oiota\, D(N,K)\ . }
$$

\smallskip

As a consequence we see that for $G := M = N$ the diagonal map $G
\longrightarrow G \times G$ induces a continuous map
$$
D(G,K)\ \mathop{\longrightarrow}\limits^{{\rm diag}_\ast}\ D(G \times
G,K) = D(G,K)\, \oiota\, D(G,K)\ .
$$
It has all the usual properties of a comultiplication; this can be
checked on Dirac distributions (where it is obvious) since they
generate dense subspaces as recalled in Remark 1.1.i.

Let now $X$ and $Y$ be two (left) $D(G,K)$-modules. Then $X \otimes_K
Y$ is a (left) $D(G,K) \otimes_K D(G,K)$-module in the obvious way. We
suppose that:

-- ${\rm dim}_K Y < \infty$, and

-- the $D(G,K)$-action on $Y$ is continuous.

Then the annihilator ideal ${\rm ann}(Y) \subseteq D(G,K)$ of $Y$ is
closed and of finite $K$-codimension. Moreover $X \otimes_K Y$ in fact
is a $D(G,K) \otimes_K D(G,K)/{\rm ann}(Y)$-module. Via the map
$$
\matrix{
D(G,K)\ \mathop{\longrightarrow}\limits^{{\rm diag}_\ast}\ D(G,K)\,
\oiota\, D(G,K)\ \mathop{\longrightarrow}\limits^{\rm pr} & D(G,K)\,
\oiota\, D(G,K)/{\rm ann}(Y) \cr
& \| \cr\cr & D(G,K) \otimes_K D(G,K)/{\rm ann}(Y) }
$$
the tensor product $X \otimes_K Y$ acquires a natural ``diagonal''
(left) $D(G,K)$-module structure.

If $X$ carries a locally convex complete Hausdorff topology with
respect to which the $D(G,K)$-action is separately continuous then the
diagonal action of $D(G,K)$ on $X\, \oiota\ Y = X\, \opi\, Y = X
\otimes_K Y$ is the unique separately continuous action which extends
the obvious diagonal action of the group $G$.

\smallskip

{\it Remark A.4:} In [ST1] we had to refer to the unpublished Diploma
thesis of F\'eaux de Lacroix for the fact that the convolution product
on $D(G,K)$ is well defined and separately continuous. But Prop. A.3
above implies in fact a more precise result. The convolution product
is induced by the multiplication map $G \times G \longrightarrow G$
via
$$
D(G,K) \times D(G,K) \rightarrow D(G,K)\,\oiota\, D(G,K) = D(G
\times G,K)\ \mathop{\longrightarrow}\limits^{\rm mult_\ast}\ D(G,K)\ .
$$

\medskip

{\bf 4. Auslander duality - derived categories}

\medskip

As a consequence of Cor. 3.7 we have the commutative diagram:
$$
\xymatrix{
  D^b(\Rep_K^{\infty,a}(G)) \ar[d]_{\widetilde{.}\,[-d]} \ar[r]^{\ .^\ast}
  & D^b_{\Cscr_G}(\Mscr_G) \ar[d]^{R\Hom_{D(G,K)}(.,\Dscr_K(G)
  \otimes_K \dfr_G)} \\
  D^b(\Rep_K^{\infty,a}(G)) \ar[r]^{\ .^\ast} & D^+(\Mscr_G)\ .
  }
$$
In this section we do always {\bf assume} that $L =
\Qdss_p$. We will see that then we can replace, as a relatively
straightforward consequence of the results in [ST2]\S8, the lower
right corner in the above diagram by $D^b_{\Cscr_G}(\Mscr_G)$.

For technical reasons we need to fix a compact open subgroup $H
\subseteq G$, and we let $D_r(H,K)$, for $r \in p^{\Qdss}$ with $1/p < r < 1$, be a
projective system of noetherian $K$-Banach algebras which exhibits
$D(H,K)$ as a Fr\'echet-Stein algebra (cf. [ST2]\S5). In the following
we will make constant use of the fact, without pointing it out each
time again, that the ring extensions $D(H,K) \longrightarrow D_r(H,K)$
are flat ([ST2] Remark 3.2). This means that the base extension
functor $D_r(H,K) \otimes_{D(H,K)} .$ is exact and hence passes
directly to the derived categories. As a piece of additional notation
we need the bounded derived category $D^b_{fg}(D_r(H,K))$ of all
complexes of (left) $D_r(H,K)$-modules whose cohomology modules are
all finitely generated and vanish in all but finitely many degrees.

\smallskip

{\bf Proposition 4.1:} {\it For any complex $X^\cdot$ in
$D^b_{\Cscr_G}(\Mscr_G)$ the $D(G,K)$-modules
$\Ext^\ast_{D(G,K)}(X^\cdot,\Dscr_K(G))$ are coadmissible and vanish
in all but finitely many degrees;

ii. for any $r$ the diagram
$$
\xymatrix{
  D^b_{\Cscr_G}(\Mscr_G) \ar[rrrr]^{R\Hom_{D(G,K)}(.,\Dscr_K(G))} \ar[d]_{\rm forget}
  & & & &
  D^b_{\Cscr_G}(\Mscr_G) \ar[d]_{\rm forget} \\
  D^b_{\Cscr_H}(\Mscr_H) \ar[rrrr]^{R\Hom_{D(H,K)}(.,D(H,K))}
  \ar[d]_{D_r(H,K) \otimes_{D(H,K)}}& & & &
  D^b_{\Cscr_H}(\Mscr_H) \ar[d]^{D_r(H,K) \otimes_{D(H,K)}} \\
  D^b_{fg}(D_r(H,K)) \ar[rrrr]^{R\Hom_{D_r(H,K)}(.,D_r(H,K))} & & & &
  D^b_{fg}(D_r(H,K))   }
$$
is commutative (up to natural isomorphism).}

Proof: i. Suppose first that $X^\cdot$ is a single coadmissible
$D(G,K)$-module $X$ concentrated in degree zero. By Prop. 2.3 we have
a $D(H,K)$-equivariant isomorphism
$$
\Ext^\ast_{D(G,K)}(X,\Dscr_K(G)) \cong \Ext^\ast_{D(H,K)}(X,D(H,K))\ .
$$
According to [ST2] Lemma 8.4 and Thm. 8.9 the right hand side is
coadmissible as a $D(H,K)$-module and vanishes for $\ast > d$. In
addition it satisfies
$$
\matrix{
D_r(H,K) \otimes_{D(H,K)} \Ext^\ast_{D(H,K)}(X,D(H,K)) \hfill\cr
\qquad\qquad\qquad\qquad \cong
\Ext^\ast_{D_r(H,K)}(D_r(H,K) \otimes_{D(H,K)} X,D_r(H,K)) \ . }
$$
Keeping in mind the fact that $\Cscr_G$ is an abelian subcategory of
$\Mscr_G$ closed under extensions the case of a general $X^\cdot$ in
$D^b_{\Cscr_G}(\Mscr_G)$ follows now by a straightforward use of the
hypercohomology spectral sequence
$$
E_2^{a,b} = \Ext^a_{D(G,K)}(h^{-b}(X^\cdot),\Dscr_K(G))
\Longrightarrow E^{a+b} = \Ext^{a+b}_{D(G,K)}(X^\cdot,\Dscr_K(G))\ .
$$
It also shows that the formula
$$
\matrix{
D_r(H,K) \otimes_{D(H,K)} \Ext^\ast_{D(H,K)}(X^\cdot,D(H,K))
\hfill\cr \qquad\qquad\qquad\qquad \cong
\Ext^\ast_{D_r(H,K)}(D_r(H,K) \otimes_{D(H,K)} X^\cdot,D_r(H,K)) \ .}
\leqno{(\ast)}
$$
holds in general.

ii. We first note that the first and second, resp. the third,
horizontal arrow is well defined by assertion i., resp. by [ST2] Thm.
8.9. The upper square is commutative up to a natural isomorphism by
Lemma 2.2 (and the fact, already observed in the proof of Prop. 2.3,
that any projective $D(G,K)$-module is projective over $D(H,K)$ as
well). The lower square is commutative by the formula $(\ast)$ above.

\smallskip

The twist functor $. \otimes_K \dfr_G : \Mscr_G \longrightarrow
\Mscr_G$ is an auto-equivalence which respects the subcategory
$\Cscr_G$ (cf. [Eme] Prop. 6.1.5).

\smallskip

{\bf Corollary 4.2:} {\it The diagram
$$
\xymatrix{
  D^b(\Rep_K^{\infty,a}(G)) \ar[d]_{\widetilde{.}\,[-d]} \ar[r]^{\ .^\ast}
  & D^b_{\Cscr_G}(\Mscr_G) \ar[d]^{R\Hom_{D(G,K)}(.,\Dscr_K(G)
  \otimes_K \dfr_G)} \\
  D^b(\Rep_K^{\infty,a}(G)) \ar[r]^{\ .^\ast} & D^b_{\Cscr_G}(\Mscr_G)
  }
$$
is commutative.}

Proof: What remains to be shown, because of Cor. 3.7, is that for any
$X^\cdot$ in $D^b_{\Cscr_G}(\Mscr_G)$ the $D(G,K)$-modules
$\Ext^\ast_{D(G,K)}(X^\cdot,\Dscr_K(G) \otimes_K \dfr_G)$ are
coadmissible and vanish in all but finitely many degrees. But we may
write $X^\cdot = Y^\cdot \otimes_K \dfr_G$ for some other complex
$Y^\cdot$ in $D^b_{\Cscr_G}(\Mscr_G)$ and then have the natural
isomorphism
$$
\Ext^\ast_{D(G,K)}(X^\cdot,\Dscr_K(G) \otimes_K \dfr_G) \cong
\Ext^\ast_{D(G,K)}(Y^\cdot,\Dscr_K(G)) \ .
$$
Since the action of $D(G,K)$ on $\dfr_G$ from the right is trivial
this isomorphism in particular is $D(G,K)$-equivariant. The claim
therefore is a consequence of Prop. 4.1.i.

\smallskip

The left perpendicular arrow in the diagram of Cor. 4.4 is an
anti-involution. What about the right perpendicular arrow?

\smallskip

{\bf Proposition 4.3:} {\it The natural transformation
$$
X^\cdot \mathop{\longrightarrow}\limits^{\cong}
R\Hom_{D(G,K)}(R\Hom_{D(G,K)}(X^\cdot,\Dscr_K(G)),\Dscr_K(G))
$$
on $D^b_{\Cscr_G}(\Mscr_G)$ is an isomorphism.}

Proof: {\it Step 1:} The natural transformation
$$
Y^\cdot \mathop{\longrightarrow}\limits^{\cong}
R\Hom_{D_r(H,K)}(R\Hom_{D_r(H,K)}(Y^\cdot,D_r(H,K)),D_r(H,K))
$$
on $D^b_{fg}(D_r(H,K))$ is an isomorphism.

Since $D_r(H,K)$ has finite global dimension $Y^\cdot$ can be
represented by a perfect complex $P^\cdot$, i.e., by a complex
consisting of finitely generated projective $D_r(H,K)$-modules at most
finitely many of which are nonzero (cf. [Ha2] III.12.3). The assertion
therefore is equivalent to the natural homomorphism of complexes
$$
P^\cdot \mathop{\longrightarrow}\limits^{\sim}
\Hom^\cdot_{D_r(H,K)}(\Hom^\cdot_{D_r(H,K)}(P^\cdot,D_r(H,K)),D_r(H,K))
$$
being a quasi-isomorphism. But this is a well known fact and can be
seen by a straightforward double complex argument.

{\it Step 2:} The natural transformation
$$
X^\cdot \mathop{\longrightarrow}\limits^{\cong}
R\Hom_{D(H,K)}(R\Hom_{D(H,K)}(X^\cdot,D(H,K)),D(H,K))
$$
on $D^b_{\Cscr_H}(\Mscr_H)$ is an isomorphism.

We fix a projective resolution $P^\cdot
\mathop{\longrightarrow}\limits^{\sim} X^\cdot$ and an injective
resolution $D(H,K) \mathop{\longrightarrow}\limits^{\sim}\break
I^\cdot$ in $\Mscr_H$. We have to show that the natural homomorphism
of complexes
$$
P^\cdot \mathop{\longrightarrow}\limits
\Hom^\cdot_{D(H,K)}(\Hom^\cdot_{D(H,K)}(P^\cdot,D(H,K)),I^\cdot)
$$
is a quasi-isomorphism, i.e., that the maps
$$
h^\ast(P^\cdot) \mathop{\longrightarrow}
h^\ast(\Hom^\cdot_{D(H,K)}(\Hom^\cdot_{D(H,K)}(P^\cdot,D(H,K)),I^\cdot))
$$
are isomorphisms. Since both sides are coadmissible it suffices to
show that the maps
$$
\matrix{
D_r(H,K) \otimes_{D(H,K)} h^\ast(P^\cdot) \hfill\cr \qquad\quad
\mathop{\longrightarrow}
D_r(H,K) \otimes_{D(H,K)}
h^\ast(\Hom^\cdot_{D(H,K)}(\Hom^\cdot_{D(H,K)}(P^\cdot,D(H,K)),I^\cdot))}
$$
are isomorphisms. The left hand side obviously is equal to
$h^\ast(D_r(H,K) \otimes_{D(H,K)} P^\cdot)$. To rewrite the right hand
side we fix an injective resolution
$$
D_r(H,K) \otimes_{D(H,K)} I^\cdot
\mathop{\longrightarrow}\limits^{\sim} J^\cdot
$$
in the category of (left) $D_r(H,K)$-modules. Applying formula
$(\ast)$ in the proof of Prop. 4.1 to the complex
$\Hom^\cdot_{D(H,K)}(P^\cdot,D(H,K))$ which,  in
$D^b_{\Cscr_H}(\Mscr_H)$, represents $R\Hom_{D(H,K)}(X^\cdot,D(H,K))$
we have the natural isomorphism
$$
\matrix{
D_r(H,K) \otimes_{D(H,K)}
h^\ast(\Hom^\cdot_{D(H,K)}(\Hom^\cdot_{D(H,K)}(P^\cdot,D(H,K)),I^\cdot))
\hfill\cr \qquad
\mathop{\longrightarrow}\limits^{\cong}
h^\ast(\Hom^\cdot_{D(H,K)}(D_r(H,K) \otimes_{D(H,K)}
\Hom^\cdot_{D(H,K)}(P^\cdot,D(H,K)),J^\cdot))\ .}
$$
using formula $(\ast)$ once more we furthermore have the
quasi-isomorphism
$$
\matrix{
D_r(H,K) \otimes_{D(H,K)} \Hom^\cdot_{D(H,K)}(P^\cdot,D(H,K))
\hfill\cr \qquad\qquad\qquad
\mathop{\longrightarrow}\limits^{\sim}
\Hom^\cdot_{D_r(H,K)}(D_r(H,K) \otimes_{D(H,K)} P^\cdot,D_r(H,K)) }
$$
which induces (cf. [Bor] I.10.5) the quasi-isomorphism
$$
\matrix{
\Hom^\cdot_{D_r(H,K)}(\Hom^\cdot_{D_r(H,K)}(D_r(H,K) \otimes_{D(H,K)}
P^\cdot,D_r(H,K)),J^\cdot) \hfill\cr \qquad\qquad
\mathop{\longrightarrow}\limits^{\sim}
\Hom^\cdot_{D_r(H,K)}(D_r(H,K) \otimes_{D(H,K)}
\Hom^\cdot_{D(H,K)}(P^\cdot,D(H,K)),J^\cdot)\ . }
$$
Our above map therefore becomes a map
$$
\matrix{
h^\ast(D_r(H,K) \otimes_{D(H,K)} P^\cdot) \hfill\cr \qquad\quad
\mathop{\longrightarrow}
h^\ast(\Hom^\cdot_{D(H,K)}(\Hom^\cdot_{D(H,K)}( D_r(H,K)
\otimes_{D(H,K)} P^\cdot,D(H,K)),J^\cdot))}
$$
which is easily checked to be the natural transformation for the
complex\break $D_r(H,K) \otimes_{D(H,K)} P^\cdot$ treated in Step 1
and hence is an isomorphism by that step.

{\it Step 3:} We fix a projective resolution $P^\cdot
\mathop{\longrightarrow}\limits^{\sim} X^\cdot$ as $D(G,K)$-modules
(and a fortiori as $D(H,K)$-modules) and injective resolutions
$\Dscr_K(G) \mathop{\longrightarrow}\limits^{\sim} J^\cdot$ as
$D(G,K)$-modules and $D(H,K) \mathop{\longrightarrow}\limits^{\sim}
I^\cdot$ as $D(H,K)$-modules. On the one hand we then have, by Lemma
2.2, the isomorphism of complexes
$$
\matrix{
{\bf l} : \Hom_{D(G,K)}(P^\cdot,\Dscr_K(G)) &
\mathop{\longrightarrow}\limits^{\cong} & \Hom_{D(H,K)}(P^\cdot,D(H,K))
\cr \hfill F & \longmapsto & \ell_{G,H}\circ F \hfill\ . }
$$
On the other hand the map $\ell_{G,H}$ extends to a map
$$
\xymatrix{
  \Dscr_K(G) \ar[d]_{\ell_{G,H}} \ar[r]^{\ \sim} & J^\cdot \ar[d]^{\ell} \\
  D(H,K) \ar[r]^{\ \ \ \sim} & I^\cdot   }
$$
of complexes. The natural transformation of our assertion is, as a map
of complexes, given as
$$
P^\cdot \longrightarrow \prod_i
\Hom_{D(G,K)}(\Hom_{D(G,K)}(P^{-i},\Dscr_K(G)),J^{\cdot +i})
$$
where $P^\cdot$ maps in the obvious way into the factor
$\Hom_{D(G,K)}(\Hom_{D(G,K)}(P^\cdot,\break \Dscr_K(G)),J^0)$ with $i
= -\cdot$ of the right hand side. A straightforward computation shows
that the diagram
$$
\xymatrix{
  P^\cdot \ar[d]_{=} \ar[r] & \prod_i
\Hom_{D(G,K)}(\Hom_{D(G,K)}(P^{-i},\Dscr_K(G)),J^{\cdot +i})
 \ar[d]^{\prod \ell\circ . \circ{\bf l}^{-1}} \\
  P^\cdot \ar[r] & \prod_i
\Hom_{D(H,K)}(\Hom_{D(H,K)}(P^{-i},D(H,K)),I^{\cdot +i})   }
$$
is commutative. By Step 2 the lower horizontal arrow is a
quasi-isomorphism. Our assertion therefore is established once we show
that the right perpendicular arrow is a quasi-isomorphism. But the
latter arrow is the composite of the maps
$$
\matrix{
\prod_i \Hom_{D(G,K)}(\Hom_{D(G,K)}(P^{-i},\Dscr_K(G)),J^{\cdot +i})
\hfill\cr\cr \qquad\qquad \mathop{\longrightarrow}\limits^{\prod\ell\circ.}
\ \prod_i \Hom_{D(H,K)}(\Hom_{D(G,K)}(P^{-i},\Dscr_K(G)),I^{\cdot +i})}
$$
and
$$
\matrix{
\prod_i \Hom_{D(H,K)}(\Hom_{D(G,K)}(P^{-i},\Dscr_K(G)),I^{\cdot +i})
\hfill\cr\cr \qquad\qquad \mathop{\longrightarrow}\limits^{\prod .\circ{\bf l}^{-1}}
\ \prod_i \Hom_{D(H,K)}(\Hom_{D(G,K)}(P^{-i},\Dscr_K(G)),I^{\cdot +i})}
$$
of which the latter even is an isomorphism of complexes. The former is
a quasi-isomorphism by Prop. 2.3 applied to the complex
$R\Hom_{D(G,K)}(X^\cdot,\Dscr_K(G))$, resp. to a projective resolution
of it (observe [Bor] I.10.5).

\smallskip

{\bf Corollary 4.4:} {\it The functor
$$
R\Hom_{D(G,K)}(.,\Dscr_K(G) \otimes_K \dfr_G) : D^b_{\Cscr_G}(\Mscr_G)
\longrightarrow D^b_{\Cscr_G}(\Mscr_G)
$$
is an anti-involution.}

Proof: From Prop. 4.3 we know that the functor
$R\Hom_{D(G,K)}(.,\Dscr_K(G))$ is an anti-involution on
$D^b_{\Cscr_G}(\Mscr_G)$. To deal with the twist by $\dfr_G$ we let
$\dfr_G : G \longrightarrow K^\times$ also denote the locally analytic
character which describes the $G$-action on the left $D(G,K)$-module
$\dfr_G$. A straightforward explicit computation shows that the map
$$
\matrix{
\Hom_{D(G,K)}(X,\Dscr_K(G)) \otimes_K \dfr_G &
\mathop{\longrightarrow}\limits^{\cong} &
\Hom_{D(G,K)}(X,\Dscr_K(G) \otimes_K \dfr_G) \cr\cr
\hfill F \otimes u & \longmapsto & [x \mapsto F(x)(\dfr_G\cdot .)
\otimes u] \hfill }
$$
is, for any $D(G,K)$-module $X$, an isomorphism of (left)
$D(G,K)$-modules. Using a projective resolution $P^\cdot
\mathop{\longrightarrow}\limits^{\sim} X^\cdot$ we now may compute
$$
\matrix{
R\Hom_{D(G,K)}(R\Hom_{D(G,K)}(X^\cdot,\Dscr_K(G)
\otimes_K \dfr_G),\Dscr_K(G) \otimes_K \dfr_G)
\hfill\cr\cr \qquad \cong
R\Hom_{D(G,K)}(\Hom^\cdot_{D(G,K)}(P^\cdot,\Dscr_K(G) \otimes_K
\dfr_G),\Dscr_K(G) \otimes_K \dfr_G)
\hfill\cr\cr \qquad \cong
R\Hom_{D(G,K)}(\Hom^\cdot_{D(G,K)}(P^\cdot,\Dscr_K(G)) \otimes_K
\dfr_G,\Dscr_K(G) \otimes_K \dfr_G)
\hfill\cr\cr \qquad \cong
R\Hom_{D(G,K)}(\Hom^\cdot_{D(G,K)}(P^\cdot,\Dscr_K(G)),\Dscr_K(G))
\hfill\cr\cr \qquad \cong
R\Hom_{D(G,K)}(R\Hom_{D(G,K)}(X^\cdot,\Dscr_K(G)),\Dscr_K(G))
\hfill\cr\cr \qquad \cong
X^\cdot\ . \hfill }
$$

\smallskip

We do know now that the outer solid arrow rectangle in the diagram
$$
\xymatrix{
  D^b(\Rep_K^{\infty,a}(G)) \ar[d]_{\widetilde{.}\,[-d]} \ar[r] &
   D^b(\Rep_K^a(G)) \ar@{-->}[d] \ar[r]^{\ .'}
  & D^b_{\Cscr_G}(\Mscr_G) \ar[d]^{R\Hom_{D(G,K)}(.,\Dscr_K(G)
  \otimes_K \dfr_G)} \\
  D^b(\Rep_K^{\infty,a}(G)) \ar[r] &
   D^b(\Rep_K^a(G)) \ar[r]^{\ .'} & D^b_{\Cscr_G}(\Mscr_G)
  }
$$is commutative with the perpendicular arrows being anti-involutions.
It remains an open question whether there is a natural perpendicular
arrow in the middle which makes each square commutative.

\medskip

{\bf 5. Auslander duality - abelian categories}

\medskip

We keep in this section the assumption that $L = \Qdss_p$. For any
coadmissible module $X$ in $\Cscr_G$ we define its {\it grade} or {\it
codimension} by
$$
j_{D(G,K)}(X) := {\rm min}\{l \geq 0 : \Ext^l_{D(G,K)}(X,\Dscr_K(G))
\neq 0\}\ .
$$
As a consequence of Prop. 2.3 we have
$$
j_{D(G,K)}(X) = j_{D(H,K)}(X)
$$
for any compact open subgroup $H \subseteq G$. Hence the present
notion of codimension coincides with the one introduced in [ST2]
p.193. From [ST2] Prop. 8.11 we then know that $X$ carries a natural
dimension filtration
$$
X = \Delta^0(X) \supseteq \Delta^1(X) \supseteq\ldots\supseteq
\Delta^{d+1}(X) = 0
$$
by $D(G,K)$-submodules having the following properties:

i. Each $\Delta^l(X)$ is a coadmissible $D(G,K)$-module;

ii. a coadmissible $D(G,K)$-submodule $X' \subseteq X$ has codimension
$\geq l$ if and only if $X' \subseteq \Delta^l(X)$;

iii. $j_{D(G,K)}(X) = {\rm sup}\{l \geq 0 : \Delta^l(X) = X\}$\ ;

iv. all nonzero coadmissible $D(G,K)$-submodules of
$\Delta^l(X)/\Delta^{l+1}(X)$ have codimension $l$.

We define $\Cscr_G^l$ to be the full subcategory in $\Cscr_G$ of all
coadmissible modules of codimension $\geq l$. By the above properties
$$
\Cscr_G = \Cscr_G^0 \supseteq \Cscr_G^1 \supseteq\ldots\supseteq
\Cscr_G^{d+1} = 0
$$
is a filtration of the abelian category $\Cscr_G$ by Serre
subcategories. We therefore may form the quotient abelian categories
$\Cscr_G^l/\Cscr_G^{l+1}$. We remark that $\Cscr^\infty_G \subseteq
\Cscr_G^d$ by [ST2] Cor. 8.13. For the purposes of our study of
duality we need the following property which is implicit in the
results in [ST2] \S8.

\smallskip

{\bf Lemma 5.1:} $j_{D(G,K)}(\Ext^l_{D(G,K)}(X,\Dscr_K(G))) \geq l$.

Proof: Let $H \subseteq G$ be a fixed compact open subgroup. By Prop.
2.3 it suffices to show that $j_{D(H,K)}(\Ext^l_{D(H,K)}(X,D(H,K)))
\geq l$. With the notations of \S4 it furthermore suffices, by [ST2]
Lemma 8.4, to see that
$$
\Ext^i_{D_r(H,K)}(\Ext^l_{D_r(H,K)}(M,D_r(H,K),D_r(H,K))) = 0
$$
for $i < l$ and any finitely generated $D_r(H,K)$-module $M$. But this
is a formal consequence of the Auslander regularity of $D_r(H,K)$
([ST2] Thm. 8.9).

\smallskip

It follows that the composed functor
$$
\Ext^l_{D(G,K)}(.,\Dscr_K(G)) : \Cscr_G^l \longrightarrow \Cscr_G^l \longrightarrow
\Cscr_G^l/\Cscr_G^{l+1}
$$
is well defined and exact. Since it is zero on the Serre  subcategory
$\Cscr_G^{l+1}$ it induces, by the universal property of quotient
categories, an exact functor $\Cscr_G^l/\Cscr_G^{l+1} \longrightarrow
\Cscr_G^l/\Cscr_G^{l+1}$.

\smallskip

{\bf Proposition 5.2:} {\it The functor
$$
\Ext^l_{D(G,K)}(.,\Dscr_K(G)) : \Cscr_G^l/\Cscr_G^{l+1} \longrightarrow
\Cscr_G^l/\Cscr_G^{l+1}
$$
is an anti-involution.}

Proof: By Prop. 4.3 we have a natural isomorphism
$$
X \cong
h^0(R\Hom_{D(G,K)}(R\Hom_{D(G,K)}(X^\cdot,\Dscr_K(G)),\Dscr_K(G))) \ .
$$
The hypercohomology spectral sequence for the right hand side
therefore is of the form
$$
E_2^{a,b} =
\Ext^a_{D(G,K)}(\Ext^{-b}_{D(G,K)}(X,\Dscr_K(G)),\Dscr_K(G))
\Longrightarrow E^{a+b}
$$
with $E^{a+b} = X$ for $a+b = 0$ and $E^{a+b}= 0$ otherwise. Let
$F^\cdot X$ denote the filtration induced by this spectral sequence on
its abutment $X$, i.e., $F^a X /F^{a+1}X = E_\infty^{a,-a}$. It
follows from Lemma 5.1 that all terms of this spectral sequence lie in
$\Cscr_G^l$. Moreover, we have $E_2^{a,b} = 0$ for $a < -b$. This
implies
$$
E_2^{a,-a} \supseteq E_3^{a,-a} \supseteq\ldots\supseteq
E_{d+2}^{a,-a} = E_\infty^{a,-a} = F^a X/F^{a+1}X
$$
for any $a \geq 0$ and with all these terms vanishing for $a < l$; in
particular, $F^l X = X$. It therefore remains to show that the
composed map
$$
X \longrightarrow X/F^{l+1}X \longrightarrow E_2^{l,-l} =
\Ext^l_{D(G,K)}(\Ext^{-l}_{D(G,K)}(X,\Dscr_K(G)),\Dscr_K(G))
$$
is an isomorphism in the quotient category $\Cscr_G^l/\Cscr_G^{l+1}$.
We have
$$
E_{i+1}^{a,-a} = {\rm ker}(E_i^{a,-a} \longrightarrow
E_i^{a+i,-a-i+1})\ .
$$
But $E_i^{a+i,-a-i+1}$, hence $E_i^{a,-a}/E_{i+1}^{a,-a}$ for $i
\geq 2$, and a fortiori $E_2^{a,-a}/E_\infty^{a,-a}$ lie in
$\Cscr_G^{a+1}$, by Lemma 5.1. It follows that $F^{l+1}X$ and
$E_2^{l,-l}/E_\infty^{l,-l}$ both lie in $\Cscr_G^{l+1}$. (We remark
that a slight refinement of this argument shows that the filtration
$F^\cdot X$ actually coincides with the dimension filtration
$\Delta^\cdot(X)$.)

\medskip

{\bf 6. The locally analytic principal series}

\medskip

In this last section we want to illustrate the preceding theory by
computing explicitly the duality functors for a series of genuinely
locally analytic representations. From now on $G$ is the group of
$\Qdss_p$-rational points of a connected reductive group over
$\Qdss_p$ (in particular $L = \Qdss_p$). We fix a parabolic subgroup
$P \subseteq G$ as well as a locally analytic character $\chi : P
\longrightarrow K^\times$. We will use the same letter
for the corresponding algebra homomorphism $\chi : D(P,K)
\longrightarrow K$. The $P$-representation, resp. the $D(P,K)$-module,
given by $\chi$ will be denoted by $K_\chi$. The locally analytic
principal series representation of $G$ corresponding to $\chi$ is
given as
$$
\matrix{
\Ind_P^G(\chi) := & \hbox{vector space of all locally analytic functions}\ f
: G \longrightarrow K\hfill\cr & \hbox{such that}\ f(gp) =
\chi^{-1}(p)f(g)\ \hbox{for any}\ g \in G\ \hbox{and}\ p \in P\hfill }
$$
with $G$ acting by left translation. By [Fea] 4.1.5 the orbit maps $g
\mapsto gf$, for $f \in \Ind_P^G(\chi)$, are locally analytic. In the
following it will be technically important to fix a maximal compact
subgroup $G_0 \subseteq G$ with the property that $G= G_0P$. We set
$P_0 := G_0 \cap P$. Then ([Fea] 4.1.4) restriction of functions
induces a $G_0$-equivariant topological isomorphism
$$
\Ind_P^G(\chi) \mathop{\rightarrow}\limits^{\cong}
\Ind_{P_0}^{G_0}(\chi)
$$
(with the right hand side having the obvious meaning). The space
$\Ind_{P_0}^{G_0}(\chi)$ is closed in $C^{an}(G_0,K)$ and hence is of
compact type with its continuous dual being a Hausdorff quotient of
$D(G_0,K)$ ([ST1] Prop. 1.2(i)). Using [ST2] Lemma 3.6 we conclude
that $\Ind_{P_0}^{G_0}(\chi)$ and $\Ind_P^G(\chi)$ are admissible
$G_0$- and $G$-representations, respectively. We now compute the dual
modules. The group $P_0$ is topologically finitely generated; we fix a
finite system $\{p_i\}_{i \in I}$ of such topological generators. The
exact sequence of admissible $G_0$-representations
$$
\matrix{
0 \longrightarrow \Ind_{P_0}^{G_0}(\chi)
\mathop{\longrightarrow}\limits^{\subseteq} C^{an}(G_0,K) &
\longrightarrow & \oplus_I\; C^{an}(G_0,K) \hfill\cr\cr
\hfill f & \longmapsto & (f(.p_i) - \chi^{-1}(p_i)f)_i \hfill }
$$
dualizes into the exact sequence of $D(G_0,K)$-modules
$$
\matrix{
\hfill \oplus_I\; D(G_0,K) & \longrightarrow &
D(G_0,K) \longrightarrow \Ind_{P_0}^{G_0}(\chi)'_b \longrightarrow 0
\cr\cr
\hfill (\lambda_i)_i & \longmapsto &
\sum_i \lambda_i\delta_{p_i} - \chi^{-1}(p_i)\lambda_i \ .\hfill }
$$
Hence
$$
\Ind_{P_0}^{G_0}(\chi)'_b = D(G_0,K)/J_{P_0}
$$
where $J_{P_0}$ denotes the left ideal of $D(G_0,K)$ generated by
$\{\delta_{p_i} - \chi^{-1}(p_i)\delta_1\}_{i \in I}$. On the other
hand a completely similar reasoning shows that the very same elements
generate the kernel of the algebra homomorphism $\chi^{-1} : D(P_0,K)
\longrightarrow K$ as a left ideal. It follows that
$$
\Ind_{P_0}^{G_0}(\chi)' = D(G_0,K) \otimes_{D(P_0,K)} K_{\chi^{-1}} \
.
$$

\smallskip

{\bf Lemma 6.1:} {\it i. The natural map
$$
D(G_0,K) \otimes_{D(P_0,K)} D(P,K)
\mathop{\longrightarrow}\limits^{\cong} D(G,K)
$$
is an isomorphism of $D(G_0,K)$-$D(P,K)$-bimodules;

ii. for any (left) $D(P,K)$-module $X$ the natural map
$$
D(G_0,K) \otimes_{D(P_0,K)} X
\mathop{\longrightarrow}\limits^{\cong}
D(G,K) \otimes_{D(P,K)} X
$$
is an isomorphism of $D(G_0,K)$-modules;

iii. there is a natural isomorphism of $D(G_0,K)$-modules
$$
{\rm Tor}_\ast^{D(P_0,K)}(D(G_0,K),X) \cong  {\rm
Tor}_\ast^{D(P,K)}(D(G,K),X)
$$
for any (left) $D(P,K)$-module $X$;

iv. in the commutative diagram
$$
\xymatrix{
   D(G_0,K) \otimes_{D(P_0,K)} K_{\chi^{-1}} \ar[d]^{\cong}
   \ar[r]^{\qquad\quad\cong}
                & \Ind_{P_0}^{G_0}(\chi)' \ar[d]^{\cong}  \\
    D(G,K) \otimes_{D(P,K)} K_{\chi^{-1}}  \ar[r]^{\qquad\quad\cong}
                & \Ind_P^G(\chi)'             }
$$
all four maps are isomorphisms.}

Proof: i.  This is clear from the decompositions $D(G,K) =
\oplus_{p \in P_0\setminus P} D(G_0,K)\delta_p$  and $D(P,K) =
\oplus_{p \in P_0\setminus P} D(P_0,K)\delta_p$. ii. This is immediate
from the first assertion. iii. Since $D(P,K)$ is free as a
$D(P_0,K)$-module a projective resolution of $X$ as a $D(P,K)$-module
can be used to compute both sides of the asserted isomorphism. Hence
the present assertion is a consequence of the previous one. iv.
(Recall that the horizontal arrows are induced by dualizing the
inclusions $\Ind_{P_0}^{G_0}(\chi) \subseteq C^{an}(G_0,K)$ and
$\Ind_P^G(\chi) \subseteq C^{an}(G,K)$, respectively, the right
perpendicular arrow is the dual of the restriction map, and the left
perpendicular arrow is induced by the inclusion $D(G_0,K) \subseteq
D(G,K)$.) We know already that the upper horizontal and the right
perpendicular arrows are isomorphisms. Hence it remains to see, for
example, that the left perpendicular arrow is bijective which is a
special case of ii.

\smallskip

Our goal is to compute the $D(G,K)$-modules
$$
\Ext^\ast_{D(G,K)}(\Ind_P^G(\chi)',\Dscr_K(G)) \cong
\Ext^\ast_{D(G,K)}(D(G,K) \otimes_{D(P,K)} K_{\chi^{-1}},
\Dscr_K(G))\ .
$$
Part of this computation can be done in greater generality for any
(left) $D(P,K)$-module $X$ which satisfies the following condition:
$$
\matrix{
\hbox{As a $D(P_0,K)$-module $X$ has a resolution $P.
\mathop{\longrightarrow}\limits^{\sim} X$ by}\hfill\cr
\hbox{finitely generated projective $D(P_0,K)$-modules.}\hfill}
\leqno{\rm (FIN)}
$$
Any such $X$ is coadmissible of course. We begin with a small
technical digression into the results of [ST2] to consider
compatibilities between $D(G_0,K)$ and its subalgebra $D(P_0,K)$.

Let $G_1$ be an open normal uniform subgroup of $G_0$.  Given an
ordered basis $h_1,\ldots, h_d$ for $G_1$ we have an explicit family
of norms $||\cdot||_r$ ($1/p<r<1$ and $r \in p^{\Qdss}$) on the
algebra $D(G_1,K)$ defining its Fr\'{e}chet-Stein structure. The
algebra $D(G_0,K)$ is a free, finite rank, left (or right)
$D(G_1,K)$-module with basis given by the Dirac distributions for
coset representatives of $G_0/G_1$, and in [ST2] Thm. 5.1 we show that
the norms $||\cdot ||_r$ extend to $D(G_0,K)$ simply by writing
elements of this algebra in a basis and taking the maximum of the
norms of the coefficients.  The Banach algebras $D_{r}(G_0,K)$, for
$r\in p^{\Qdss}$ and $1/p<r<1$, are the completions of $D(G_0,K)$ with
respect to the norms $||\cdot||_r$.

Naturally, all of the constructions described here may also be carried
out for the group $P_0$. The following result shows that they may be
done ``simultaneously'' for $P_0$ and $G_0$.

\smallskip

{\bf Proposition 6.2:} {\it One can choose the family of norms
$||\cdot||_r$ on $D(G_0,K)$ so that they define the Fr\'{e}chet-Stein
structure on both $D(G_0,K)$ and its subalgebra $D(P_0,K)$, and so
that, for each $1/p<r<1$ in $p^{\Qdss}$, the completion $D_r(G_0,K)$
is flat as a left as well as a right $D_r(P_0,K)$-module.}

Proof: The result holds  more generally for an arbitrary compact
$p$-adic Lie group $H$ and a closed subgroup $Q \subseteq H$, and we
will prove it in this situation. By [DDMS] Ex. 4.14, we may find an
open normal uniform subgroup $H'$ of $H$ such that $H'\cap Q=Q'$ is
open normal uniform in $Q$. Furthermore, we may arrange that $(H')^p
\cap Q' = Q'^p$. In this situation we find an ordered basis of
topological generators $a_{1},\ldots,a_{d}$ for $H'$ such that the
first $m := \dim(Q)$ members are an ordered basis for $Q'$ (cf. [DDMS]
Prop. 1.9(iii) and Lemma 3.4). Therefore inside the explicit
realization of $D(H',K)$ as a noncommutative power series ring in the
variables $b_{i} := a_{i}-1$ for $i=1,\ldots, d$, as in [ST2] \S4, the
subalgebra $D(Q',K)$ consists of those power series involving only the
first $m$ variables.  The norms $||\cdot||_r$ for $D(H',K)$ restrict
to the corresponding ones for $D(Q',K)$. The statement and proof of
Thm. 4.5 of [ST2] shows that, on the level of graded rings, the map
from $D_{r}(Q',K)$ to $D_{r}(H',K)$ is just the inclusion of a
polynomial ring  in $m+1$ variables into one of $d+1$ variables, which
is clearly flat.  The conditions of Prop. 1.2 of [ST2] being
satisfied, this tells us that $D_{r}(Q',K)\to D_{r}(H',K)$ is also
flat, both on the left and right.

Now following the proof of [ST2] Thm. 5.1, we extend the norms on
$D(H',K)$ to $D(H,K)$ by choosing coset representatives for $H'$ in
$H$, writing
$$
D(H,K)=\oplus_{h\in H/H'} D(H',K)\delta_{h}
$$
as a free (left) $D(H',K)$ module, and, for each $r$,  taking the
maximum of the $||\cdot||_r$--norm of the components in this
representation. If  we further require that representatives in the
image of the inclusion map $Q/Q' \hookrightarrow H/H'$ be chosen from
$Q$, then the restriction of $||\cdot||_r$ to $D(Q,K)=\oplus_{h\in
Q/Q'}D(Q',K)\delta_{h}$ respects this direct sum. By the flatness
result in the uniform case, we know that $D_{r}(Q,K)\otimes_{
D_{r}(Q',K)} D_{r}(H',K)$ is flat as a left $D_{r}(Q,K)$-module.
However, a computation with coset representatives shows that
$D_{r}(H,K)$ is a finite direct sum of copies of
$D_r(Q,K)\otimes_{D_{r}(Q',K)}D_{r}(H',K)$, and is therefore flat as a
left $D_{r}(Q,K)$-module. The right flatness follows in the same way.

\smallskip

{\bf Lemma 6.3:} {\it Suppose that the (left) $D(P,K)$-module $X$
satisfies (FIN); we then have:

i. ${\rm Tor}_\ast^{D(P,K)}(D(G,K),X) = {\rm
Tor}_\ast^{D(P_0,K)}(D(G_0,K),X) = 0$ for $\ast > 0$;

ii. there is a natural isomorphism of (right) $D(G,K)$-modules}
$$
\Ext^\ast_{D(G,K)}(D(G,K) \otimes_{D(P,K)} X,
\Dscr_K(G)) \cong \Ext^\ast_{D(P,K)}(X,\Dscr_K(G)) \ .
$$
Proof: i. Because of Lemma 6.1.iii we only need to show the vanishing
of the second Tor-groups. Using the projective resolution $P.
\mathop{\longrightarrow}\limits^{\sim} X$ from (FIN) we have
$$
{\rm Tor}_\ast^{D(P_0,K)}(D(G_0,K),X) = h_\ast(D(G_0,K)
\otimes_{D(P_0,K)} P.) \ .
$$
But $D(G_0,K) \otimes_{D(P_0,K)} P.$ is a complex of coadmissible
$D(G_0,K)$-modules. So its homology is coadmissible which means that
its vanishing can be tested on the corresponding coherent sheaves (cf.
[ST2] Cor. 3.1). By [ST2] Remark 3.2 we have
$$
D_r(G_0,K) \otimes_{D(G_0,K)} h_\ast(D(G_0,K)
\otimes_{D(P_0,K)} P.) = h_\ast(D_r(G_0,K)
\otimes_{D(P_0,K)} P.) \ .
$$
Since both ring homomorphisms $D(P_0,K) \longrightarrow D_r(P_0,K)
\longrightarrow D_r(G_0,K)$ are flat, the first one by [ST2] Remark
3.2 and the second one by Prop. 6.2, the groups $h_\ast(D_r(G_0,K)
\otimes_{D(P_0,K)} P.)$ vanish in degrees $\ast > 0$.

ii. The assertion is obvious in degree $\ast = 0$. The general case
follows from this if we use a projective resolution $Q.
\mathop{\longrightarrow}\limits^{\sim} X$ of $D(P,K)$-modules to
compute the right hand side because $D(G,K) \otimes_{D(P,K)} Q.
\mathop{\longrightarrow}\limits^{\sim} D(G,K) \otimes_{D(P,K)} X$
then, by i., is a projective resolution which we may use to compute
the left hand side.

\smallskip

{\bf Proposition 6.4:} {\it Suppose that the (left) $D(P,K)$-module
$X$ as well as all the (right) $D(P,K)$-modules
$\Ext^\ast_{D(P,K)}(X,\Dscr_K(P))$ satisfy (FIN); then there is a
natural isomorphism of (right) $D(G,K)$-modules}
$$
\matrix{
\Ext^\ast_{D(G,K)}(D(G,K) \otimes_{D(P,K)} X,
\Dscr_K(G)) \cong \hfill\cr
\qquad\qquad\qquad\qquad\qquad\qquad\qquad
\Ext^\ast_{D(P,K)}(X,\Dscr_K(P)) \otimes_{D(P,K)} D(G,K) \ .}
$$
Proof: Since $\Dscr_K(P) \subseteq \Dscr_K(G)$ we have, by
functoriality, the homomorphism of right $D(G,K)$-modules
$$
\matrix{
\Ext^\ast_{D(P,K)}(X,\Dscr_K(P)) \otimes_{D(P,K)} D(G,K) & \longrightarrow
& \Ext^\ast_{D(P,K)}(X,\Dscr_K(G)) \hfill\cr
\hfill e \otimes \lambda & \longmapsto &
\Ext^\ast_{D(P,K)}(X,.\cdot\lambda)(e)\ . \hfill}
$$
We claim that this map in fact is an isomorphism, which by Lemma
6.3.ii suffices for our assertion. For this we consider the diagram:
$$
\xymatrix{
  \Ext^\ast_{D(P,K)}(X,\Dscr_K(P)) \otimes_{D(P,K)} D(G,K)
     \ar[r] & \Ext^\ast_{D(P,K)}(X,\Dscr_K(G)) \ar[dd]_{\cong}^{\Ext^\ast(X,\ell_{G,G_0})} \\
  \Ext^\ast_{D(P,K)}(X,\Dscr_K(P)) \otimes_{D(P_0,K)} D(G_0,K)
     \ar[u]_{\cong} \ar[d]^{\cong}_{\Ext^\ast(X,\ell_{P,P_0})\otimes id} &  \\
  \Ext^\ast_{D(P_0,K)}(X,D(P_0,K)) \otimes_{D(P_0,K)} D(G_0,K) \ar[r] &
  \Ext^\ast_{D(P_0,K)}(X,D(G_0,K))   }
$$
The upper left perpendicular arrow is an isomorphism by the right
module version of Lemma 6.1.ii. The lower left and the right
perpendicular arrows are isomorphisms by Prop. 2.3, resp. an argument
entirely analogous to the proof of Prop. 2.3. The commutativity of
this diagram reduces to the commutativity of
$$
\xymatrix{
  \Dscr_K(P) \ar[d]_{\ell_{P,P_0}} \ar[r]^{\cdot \mu} & \Dscr_K(G) \ar[d]^{\ell_{G,G_0}} \\
  D(P_0,K) \ar[r]^{\cdot \mu} & D(G_0,K)   }
$$
for any $\mu \in D(G_0,K)$ which is obvious. Our claim therefore
reduces to the lower horizontal arrow in the first diagram being an
isomorphism. Let $P. \mathop{\longrightarrow}\limits^{\sim} X$ be the
projective resolution from (FIN). We have to show that
$$
\xymatrix{
 h^\ast(\Hom_{D(P_0,K)}(P.,D(P_0,K))) \otimes_{D(P_0,K)} D(G_0,K)
     \ar[d] \\
 h^\ast(\Hom_{D(P_0,K)}(P.,D(G_0,K))) }
$$
is an isomorphism. Since each term in the complex $P.$ is finitely
generated projective we have
$$
\Hom_{D(P_0,K)}(P.,D(G_0,K)) = \Hom_{D(P_0,K)}(P.,D(P_0,K))
\otimes_{D(P_0,K)} D(G_0,K) \ .
$$
This finally reduces our claim to the statement that the natural map
$$
h^\ast(Y^\cdot) \otimes_{D(P_0,K)} D(G_0,K) \longrightarrow
h^\ast(Y^\cdot \otimes_{D(P_0,K)} D(G_0,K))
$$
is an isomorphism for the complex $Y^\cdot :=
\Hom_{D(P_0,K)}(P.,D(P_0,K))$. The complex $Y^\cdot$
consists of finitely generated projective right $D(P_0,K)$-modules. By
Prop. 4.1 it is bounded. Moreover, by assumption, its cohomology are
(right) $D(P_0,K)$-modules which satisfy (FIN). Therefore this
statement is a formal consequence, by a hypercohomology spectral
sequence argument, of the right module version of Lemma 6.3.i.

\smallskip

{\bf Proposition 6.5:} {\it i.
$\Ext^\ast_{D(G,K)}(\Ind_P^G(\chi)',\Dscr_K(G)) = 0$ for $\ast \neq
{\rm dim}(P)$;

ii. $\Ext^{{\rm dim}(P)}_{D(G,K)}(\Ind_P^G(\chi)',\Dscr_K(G)) \cong
\Ind_P^G((\chi\dfr_P)^{-1})'$.}

Proof: (We emphasize that in ii. we consider both sides, as usual, as
left $D(G,K)$-modules.) First of all we have to verify that the
$D(P,K)$-module $X := K_{\chi^{-1}}$ satisfies the assumptions of
Prop. 6.4. Let $D^c(P_0,K)$ denote the convolution algebra of
$K$-valued continuous distributions (= measures) on the compact group
$P_0$. It is the continuous dual of the Banach space of all $K$-valued
continuous functions on $P_0$. Alternatively it can be constructed as
follows. Let $o$ denote the ring of integers in $K$ and form the
completed group ring
$$
o[[P_0]] := \mathop{\rm lim}\limits_{\mathop{\longleftarrow}\limits_N}
o[P_0/N]
$$
where $N$ runs over all open normal subgroups of $P_0$. Then
$$
D^c(P_0,K) \cong K \otimes_o o[[P_0]] \ .
$$
The ring $o[[P_0]]$ is known to be noetherian by a straightforward
generalization of [Laz] V.2.2.4; hence $D^c(P_0,K)$ is noetherian. The
continuous homomorphism $\chi^{-1} : P_0 \longrightarrow o^\times
\subseteq K^\times$ extends to an algebra homomorphism $\chi^{-1} :
D^c(P_0,K) \longrightarrow K$ which allows us to view $K_{\chi^{-1}}$
as a (left) $D^c(P_0,K)$-module. Since $D^c(P_0,K)$ is noetherian we
find a resolution $Q. \mathop{\longrightarrow}\limits^{\sim}
K_{\chi^{-1}}$ by finitely generated projective $D^c(P_0,K)$-modules.
In [ST2] Thm. 5.2 we have shown that the natural ring homomorphism
$D^c(P_0,K) \longrightarrow D(P_0,K)$ is flat. For the sake of clarity
we point out that the statement of loc. cit. only says that the map
$D^c(P_0,K') \longrightarrow D(P_0,K)$ is flat for any subfield $K'
\subseteq K$ which is finite over $\Qdss_p$. But noticing that
$o[[P_0]]$ is a linearly compact $o$-module the same proof actually
gives this slightly more general result needed here. (Warning: The
ring denoted by $K[[P_0]]$ in loc. cit. in general is dense in but not
equal to $D^c(P_0,K)$.) It follows that
$$
P. := D(P_0,K) \otimes_{D^c(P_0,K)} Q.
\mathop{\longrightarrow}\limits^{\sim}
D(P_0,K) \otimes_{D^c(P_0,K)} K_{\chi^{-1}}
$$
is a resolution by finitely generated projective $D(P_0,K)$-modules.
Since, as discussed at the beginning of this section, the kernel of
$\chi^{-1} : D(P_0,K) \longrightarrow K$, resp. of $\chi^{-1} :
D^c(P_0,K) \longrightarrow K$, is the left ideal generated by the
elements $\delta_p - \chi^{-1}(p)\delta_1$, for $p \in P_0$, in the
respective ring we have
$$
D(P_0,K) \otimes_{D^c(P_0,K)} K_{\chi^{-1}} = K_{\chi^{-1}} \ .
$$
Hence $P. \mathop{\longrightarrow}\limits^{\sim} K_{\chi^{-1}}$ is a
resolution as required in (FIN). Furthermore, by the computation in
the proof of Cor. 4.4 we have right $D(P,K)$-module isomorphisms
$$
\eqalign{
\Ext^\ast_{D(P,K)}(K_{\chi^{-1}},\Dscr_K(P))
 & = \Ext^\ast_{D(P,K)}(K,\Dscr_K(P) \otimes_K K_{\chi})\cr
 & = \Ext^\ast_{D(P,K)}(K,\Dscr_K(P) \otimes_K \dfr_P) \otimes_K K_{(\chi\dfr_P)^{-1}} \ .}
$$
Applying Cor. 3.7 to the trivial (and hence smooth) $P$-representation
we obtain
$$
\Ext^\ast_{D(P,K)}(K,\Dscr_K(P) \otimes_K \dfr_P) =
\left\{\matrix{K \hfill &
\hbox{if}\ \ast = {\rm dim}(P),\hfill\cr
0 \hfill & \hbox{otherwise.} \hfill\cr}\right.
$$
as right $D(P,K)$-modules. We see that
$$
\Ext^\ast_{D(P,K)}(K_{\chi^{-1}},\Dscr_K(P)) =
\left\{\matrix{K_{(\chi\dfr_P)^{-1}} \hfill &
\hbox{if}\ \ast = {\rm dim}(P),\hfill\cr
0 \hfill & \hbox{otherwise.} \hfill\cr}\right.
$$
as right $D(P,K)$-modules. In particular these $D(P,K)$-modules
satisfy (FIN) as well. We therefore may apply Prop. 6.4 and, together
with Lemma 6.1.iv, we obtain
$$
\Ext^\ast_{D(G,K)}(\Ind_P^G(\chi)',\Dscr_K(G)) \cong
\Ext^\ast_{D(P,K)}(K_{\chi^{-1}},\Dscr_K(P)) \otimes_{D(P,K)} D(G,K)
$$
as right $D(G,K)$-modules. Inserting into this the previous
computation, converting right into left modules, and using Lemma
6.1.iv once more establishes the assertion.

\bigskip

{\bf References}

\parindent=23truept

\ref{[Bor]} Borel A. et al.: Algebraic $D$-Modules. Orlando: Academic
Press 1987

\ref{[BW]} Borel A., Wallach N.: Continuous cohomology, Discrete Subgroups,
and Representations of Reductive Groups. Ann. Math. Studies 94.
Princeton Univ. Press 1980

\ref{[B-GAL]} Bourbaki, N.: Groupes et alg\`ebres de Lie, Chap. 1-3. Paris:
Hermann 1971, 1972

\ref{[CE]} Cartan H., Eilenberg S.: Homological Algebra. Princeton
Univ. Press 1956

\ref{[Cas]} Casselman W.: Introduction to the theory of admissible
representations of $\wp$-adic reductive groups. Preprint

\ref{[Dix]} Dixmier J.: Enveloping Algebras. Amsterdam: North Holland
1977

\ref{[DDMS]} Dixon J.D., du Sautoy M.P.F., Mann A., Segal D.: Analytic
Pro-$p$-Groups (2nd Edition). Cambridge Univ. Press 1999

\ref{[Eme]} Emerton M.: Locally analytic vectors in representations of
non-archi- medean locally $p$-adic analytic groups. Preprint 2003

\ref{[Fea]} F\'eaux de Lacroix C. T.: Einige Resultate \"uber die topologischen
Dar- stellungen $p$-adischer Liegruppen auf unendlich dimensionalen
Vektor- r\"aumen \"uber einem $p$-adischen K\"orper. Thesis, K\"oln
1997, Schriftenreihe Math. Inst. Univ. M\"unster, 3. Serie, Heft 23,
pp. 1-111 (1999)

\ref{[Gro]} Grothendieck A.: Produit tensoriels topologiques et
espaces nucl\'eaires. Memoirs AMS 16 (1966)

\ref{[Har]} Hartshorne R.: Residues and Duality. Lect. Notes Math. 20.
Berlin-Heidelberg-New York: Springer 1966

\ref{[Ha2]} Hartshorne R.: Algebraic Geometry. Berlin-Heidelberg-New York:
\break Springer 1977

\ref{[Kem]} Kempf G. R.: The Ext-dual of a Verma module is a Verma module.
J. Pure Appl. Algebra 75, no. 1, 47--49 (1991)

\ref{[Lam]} Lam T. H.: Lectures on Modules and Rings. Berlin-Heidelberg-New York:
Springer 1999

\ref{[Laz]} Lazard M.: Groupes analytiques $p$-adique. Publ. Math.
IHES 26, 389-603 (1965)

\ref{[Sch]} Schneider P.: $p$-adische Analysis. Course at M\"unster in
2000.\hfill\break
http://wwwmath.uni-muenster.de/math/u/schneider/publ/lectnotes

\ref{[NFA]} Schneider P.: Nonarchimedean Functional Analysis. Berlin-Heidelberg-New York:
Springer 2001

\ref{[ST1]} Schneider P., Teitelbaum J.: Locally analytic
distributions and $p$-adic representation theory, with applications to
$GL_2$. J. AMS 15, 443-468 (2002)

\ref{[ST2]} Schneider P., Teitelbaum J.: Algebras of $p$-adic
distributions and admissible representations. Invent. math. 153,
145-196 (2003)

\bigskip

\parindent=0pt

Peter Schneider\hfill\break Mathematisches Institut\hfill\break
Westf\"alische Wilhelms-Universit\"at M\"unster\hfill\break
Einsteinstr. 62\hfill\break D-48149 M\"unster, Germany\hfill\break
pschnei@math.uni-muenster.de\hfill\break
http://www.uni-muenster.de/math/u/schneider\hfill

\noindent
Jeremy Teitelbaum\hfill\break Department of Mathematics, Statistics,
and Computer Science (M/C 249)\hfill\break University of Illinois at
Chicago\hfill\break 851 S. Morgan St.\hfill\break Chicago, IL 60607,
USA\hfill\break jeremy@uic.edu\hfill\break
http://raphael.math.uic.edu/$\sim$jeremy\hfill

\end